\documentclass[12pt]{amsart}
\usepackage{amsbsy,amssymb,amscd,amsfonts,latexsym,amstext,delarray,
amsmath}
\newtheorem{thm}{Theorem}[section] 
\newtheorem{prop}[thm]{Proposition}       
\newtheorem{cor}[thm]{Corollary}           
\newtheorem{lem}[thm]{Lemma}                   

\newtheorem{defn}{Definition}[section]
\newtheorem{rem}[thm]{Remark}

\numberwithin{equation}{section}

\def\R{{\mathbb R}} 
\def\A{{\bf A}} 
\def\Pr{{\bf P}}

\def\Ql{{\mathbb Q}_{\ell}}

\def\Q{{\mathbb Q}}
\def\C{{\mathbb C}}  
\def\Z{{\mathbb Z}}
\def\F{{\mathbb F}}       
\def\mod{{\rm mod}}    
\def\gcd{{\rm gcd}}    
\def\Aut{{\rm Aut}}    
\def\Gal{{\rm Gal}}    
\def\fr{{\rm Fr}}      
\def\id{{\rm Id}}      
\def\tr{{\rm Trace}}   
\def\Ind{{\rm Ind}}    
\def\Re{{\rm Re}}      

\newcommand{\ie}{{\it i.e.\/}\ }
\newcommand{\eg}{{\it e.g.\/}\ }
\newcommand{\cf}{{\it cf.\/}\ }
\newcommand{\resp}{{\it resp.\/}\ }
\newcommand{\op}{{\it op.cit.\/}\ }

\begin{document}

\title[Geometry and arithmetic on a quintic threefold]{Geometry and arithmetic on a quintic threefold}

\author[Caterina Consani]{Caterina Consani$^{\text{\dag}}$}
\author[Jasper Scholten]{Jasper Scholten}
\thanks{\noindent$^{\text{\dag}}$Partially supported by the NSF grant DMS-9701302 and by the NSERC grant 72016789}
\address{C.~Consani: Department of Mathematics \\
University of Toronto \\
Toronto, ON.\ \ M5S 3G3 CA
}
\email{kc@math.toronto.edu}
\address{J.~Scholten:  Department of mathematics\\
Katholieke Universiteit Nijmegen\\
Postbus 9010, 6500 GL, Nijmegen,\ \ the Netherlands}

\email{scholten@sci.kun.nl}
\maketitle

\section*{Introduction}\label{00}

In this paper we investigate the geometry and the arithmetic of a non-rigid Calabi-Yau threefold $\overline X$ with $120$ nodes.
This variety is the projective closure in $\Pr^4_{\C}$ of a special fiber 
in a family of affine threefolds. The family is obtained as a self-fiber product of a fibration of affine curves defined by generalized Chebyshev polynomials in two variables. 
These polynomials are particular kind of symmetric functions with the remarkable property of taking only three different values at their critical points. They naturally 
generalize in higher dimensions  the classical Chebyshev polynomials in one variable. The self-fiber
product of a one parameter family of degree five generalized Chebyshev polynomials $P_5(x_1,x_2) = t$ is an affine fibration endowed with an action of the dihedral group of order
eight. The fiber at the origin $X$ has only non-degenerate singular points as singularities. The 
presence of a group action on the 
fibers of the self-fiber product (and in particular on $X$) makes possible the computation of the Hodge numbers of  a desingularization $\tilde X$
of $\overline X$. Computations show that  $\tilde X$ is a Calabi-Yau threefold with $h^2(\tilde X) = h^{1,1}(\tilde X) = 141$ and $h^3(\tilde X) = 4$ (\cf corollaries~\ref{cor} and \ref{dimensions}). 

The particular position
of the singular points on $X$ --that is a consequence of the symmetries of the Chebyshev polynomials-- allows the presence of a divisor
on the threefold through the nodes and not homologous to a multiple of a generic hyperplane section. 
This cycle is the responsible for the non-vanishing of $h^{2,1}(\tilde X)$.

The affine threefold $X$ has a covering $Y$ whose definition is given in terms of Dickson polynomials.
More precisely, $X$ is the variety of the orbits for the action of (two copies of) the symmetric
group in three variables on $Y$. Although the geometry of $Y$ is much more complicated than that of $X$
(\eg its singular locus consists of $8750$ points and unlikely $X$, $Y$ acquires further singularities at infinity),
the presence of this covering plays a relevant role in the proof of the existence of  correspondences on $X$. 
These correspondences  provide a
reasonable explanation for some of the arithmetic results we prove in the paper.

Over the field of the rational numbers, $\tilde X$ has good reduction outside the set $\{2,3,5\}$ (\cf section~\ref{J1}).  One of the main goals of this paper is the investigation of the Galois representation 
$\rho: {\rm Gal}(\overline\Q/\Q) \to{\rm GL}(H^3(\tilde X_{\overline\Q}, {\Q}_{\ell}))$ and the L-function ${\rm L}(H^3(\tilde X_{\overline\Q}, {\Q}_{\ell}),s)$ ($\ell$ a rational prime) of the resolved variety $\tilde X$. 
 This study is rather complicated and it is  
accomplished in few steps. First, we compute the local Euler factors at many (rational) primes of good reduction (\cf paragraph~\ref{J2}). This is done by counting the points of $\tilde X$ over various finite fields and using the 
Lefschetz trace formula. The shape of these polynomials and their splitting behavior in the real quadratic extension
 $F = \Q(\sqrt 5)$,  clearly suggest that the scalars extension $\rho \otimes \Q_\ell(\sqrt 5)$ is induced  by a two--dimensional representation $\sigma$ of the Galois group ${\rm Gal}(\overline\Q/F)$. We prove this in paragraph \ref{J3} (\cf~theorem~\ref{Jaspth}). 
The existence of  algebraic correspondences on $X$ defined over $F$ (\cf the discussion at the end of paragraph \ref{J3}), provides  a reason for the reducibility of the restriction $(\rho \otimes \overline\Q_\ell)|_{{\rm Gal}(\overline\Q/F)}$  and another explanation for the induceness.

Following the theory of Jacquet and Langlands one would predict the modularity of $\sigma$, namely the existence of a Hilbert modular form $\mathfrak f$ over $F$ whose associated Galois representation $\sigma_{\mathfrak f}$ is equivalent to $\sigma$. We deduce the expected weight ${\bf k} = (2,4)$ and the level $\mathfrak N = 5$ of $\mathfrak f$  from the shape of the (conjectural) functional equation of the L--function of the threefold (\cf paragraph~\ref{J4}). This procedure is possible because one knows the invariance of the L--function for induced representations. 
With the help of a computer, we tested numerically the functional equation, under the assumption of its analytic continuation. 

The existence of a Hilbert cusp form with the correct behavior is deduced from the theorem of Jacquet and Langlands. The properties of $\mathfrak f$ are characterized on the automorphic side, by those of the corresponding automorphic representation $\pi_{\mathfrak f}$ of $G_B(\A_F)$: the adeles of the algebraic group $G_B$ associated to (the multiplicative group of) a definite quaternion algebra $B$ over $F$. The definition of this algebra and the description of a correlated Eichler order are derived from the expected description of the conductor in the functional equation (\cf section~\ref{3}).

The action of the Hecke operators on the space of the quaternionic forms is given by the Brandt matrices. These square matrices are naturally associated to an Eichler order and they describe the representation of the Hecke operators acting on a space of theta series attached to the norm form of the algebra $B$.
In our case, it is known (\cf~\cite{W}) that these series provide 
a basis for the space of Hilbert modular forms to which  $\mathfrak f$ 
belongs. In the paper we provide an explicit description of the algorithm  
used for the definition of the Brandt matrices (\cf~paragraph~\ref{6}). 
It is known that these matrices generate a commutative semisimple 
ring and they satisfy the same identities fulfilled by the classical Hecke 
operators. Hence, through a process of simultaneus diagonalization of them 
we provide a proof of the existence of a Hilbert modular form $\mathfrak f'$ 
with properties analogous to those of $\mathfrak f$ (\cf paragraph~\ref{6}). 
The function $\mathfrak f'$ is represented by a common eigenvector 
of the Brandt matrices and its eigenvalues are known to be equal to the traces 
of Frobenius of the associated Galois representation $\sigma_{\mathfrak f'}$.

In our case, all the computed eigenvalues are in fact equal to
the traces of Frobenius of the geometric Galois representation $\sigma$.
This suggests that $\sigma$ and $\sigma_{\mathfrak f'}$ are isomorphic upto
semisimplification.  By a theorem of 
Faltings (\cf~\cite{F}), two semisimple
Galois representations are isomorphic if their traces of Frobenius are 
equal for sufficiently many primes within a finite set. 
Unfortunately,  we are not able to make this test
set small enough to be contained in the set of primes for which we
can compute the traces. A method described in \cite{RL} allows one
to make the test set very small if one knows that the traces of
Frobenius are even. In our case all traces appear to be even, but we can only
prove that for the geometric Galois representation $\sigma$.



There are at least two further problems that we have not investigated in this paper but whose study seems natural to pursue. 

The very rich geometry of $\tilde X$ makes this variety a good example for testing some of the open conjectures in arithmetic and in the study of the algebraic cycles. The theory of Jacquet and Langlands predicts the modularity of the  4-dimensional selfdual Galois representation $\rho$. More precisely, we expect to support evidence for the existence of a Siegel modular form of weight 3 and genus 2, defined on some subgroup $\Gamma_0(N)$ of ${\rm Sp}(4,\Z)$, whose associated Galois representation is equivalent to $\rho$. It is concievable, although we have not checked this in details, that this modular form arises as (theta) lifting of the Hilbert modular form $\mathfrak f'$.

 A second theme of investigation is the verification of Bloch's conjecture that relates in our case, the order of zero of the L--function at the center of the critical strip and the nonvanishing of the Griffiths group of the algebraic cycles on $\tilde X$ homologous to zero modulo algebraic equivalence. From the numerical analysis made on the functional equation, we were able to test the vanishing of the L--function ${\rm L}(H^3,s)$ at $s = 2$, exactly on the first order. In accord with Bloch's conjecture, this first order of zero would be explained geometrically by the presence of a codimension--two algebraic cycle on $\tilde X$, whose class in the Griffiths group is not trivial. We expect that the construction described by Kimura in \cite{K} may be helpful for the definition of this cycle. \vspace{.1in}

In the following, we give a brief account on the arguments illustrated in the various paragraphs.

In the first two sections we define the threefold fibration in which $X$ appears as a special fiber and we compare our construction with others  made in \cite{vGW} and \cite{vGW1}.
In the third paragraph we describe the group action on $X$ using which we compute the Hodge numbers.
The sections~\ref{J1},\ref{J2},\ref{J3} and \ref{J4} include the study of the L-function and the Galois representation of $\tilde X$. 
The paragraphs~\ref{3},\ref{4} are dedicated  to the definition and the properties of the quaternion  algebra $B$ and the Eichler order. The final sections~\ref{5},\ref{6} report on the description of the algorithm  used for the definition of the Brandt matrices and for their simultaneus diagonalization. 

In the paper we have included few tables which resume the high number of computations we have made. In many situations, the help of an appropriate computer program (Mathematica, Maple, Pari, Magma) was fundamental.  
\subsection*{\it Acknowledgments} 

This paper is the result of a rather long period of research  dedicated to the study of the Langlands program. Our intention was to approach this complicated theory via the study of a hypersurface whose  geometry and arithmetic turned out to be very rich. In many occasions, we referred to the experienced advice of different mathematicians to obtain answers or for explaining our computations and  listen their comments. This is a heartful thank to all who have been contacted by us either in person or by e-mail. Among these people we would like to mention Juliusz Brzezinski and Elise Bjrkholdt to whom we often referred for questions regarding the theory of orders in quaternion algebras. We thank Fred Diamond who  explained to us the approach to the study of the Hilbert modular forms via functions defined on a quaternion algebra. Many thanks are due to Dinakar Ramakrishnan for his advice and because he suggested the lecture of the thesis \cite{So}. This paper became for us an inestimable source of reference, especially in the process of setting up the algorithm for the computation of the Brandt matrices.  Our thanks go to Spencer Bloch and Jaap Top for the support and the encouragement we received from them in these years and to Jim Cogdell, Pierre Deligne and Bert van Geemen who, in many occasions, patiently listened our ideas and replied with interesting remarks. The first author gratefully thanks the Max-Planck Institut f\"ur Mathematik and the Institut des Hautes \'Etudes Scientifiques for the kind hospitality received during her stay there, respectively in June 1999 and January 2000. She also acknowledges the support received by the Institute for Advanced Study in the year 1999-2000.
The second author wishes to thank the first author for asking him to join
this project, and for the hospitality received during various visits.

\section{A skew-pentagon configuration of lines in the affine
plane.}\label{0} 

In this paragraph we review the description of a family of
affine, plane, quintic curves introduced by van Geemen and Werner (\cf~\cite{vGW} 
and \cite{vGW1}) in relation with the study of a particular
configuration of lines in the plane.
The authors called this configuration a
skew pentagon (five lines meeting only in pairs and stable under the involution 
$(x,y) \mapsto (x,-y)$) as it generalizes  the construction of the regular pentagon defined by
Hirzebruch in \cite{H}. The equation of a skew pentagon depends upon two parameters $a,~b$ as follows
\[
F_{a,b}(x,y) = (x+a)(y^2-x^2)(y^2-b(x+1)^2).
\]
It is easy to verify that $F_{a,b}$ has 10 critical points (\ie all of its partial derivatives vanish)
where the five lines
mutually intersect. The polynomial $F_{a,b}$ has 2 further critical
points on the x-axis. If one adds the condition that
\begin{equation}\label{condition}
\text{the value of $F_{a,b}$ at the remaining 4 
critical points (not on the lines) is the same,
}
\end{equation}
it turns out that the parameters $a,b$ must satisfy the  
quartic equation 
\begin{equation}\label{quartic}
a^2(b-1)^2-2ab(b-1)+b(b-5) = 0.
\end{equation}
Van Geemen and Werner noticed that the solutions of \eqref{quartic} are
parametrized via the introduction of an additional parameter $t$
\[
a = \frac{t(t+5)}{t^2-5},\quad b = \frac{t^2}{5}.
\]

Set $H_t(x,y) := F_{a(t),b(t)}(x,y)$. Then, $H_t(x,y)$ defines a
family of  plane curves parametrized by the rational quartic
\eqref{quartic}. Each element of this family
defines a skew pentagon configuration of lines in the affine plane under
the condition \eqref{condition}.

It is easy to check that for a generic value of $t$ (\ie for all $t$
except finitely many of them), the quintic   
threefold $V \subset \Pr^4$ defined as the closure of the affine variety in $\A^4$
\begin{equation}\label{Bertfamily}
H_t(x,y) - H_t(u,v) = 0
\end{equation}
has 118 nodes (rational double points): \ie $10\cdot 10 + 4 \cdot 4 + 1 \cdot 1 + 1 \cdot 1$. These points are the
only kind of singularities of $V$ and they are all
located on an affine chart (\ie the affine threefold defined by the equation \eqref{Bertfamily} does not
acquire further singular points in its closure). The condition \eqref{condition} constrains
the critical values taken by the function $H_t(x,y)$ (for a generic choice
of $t$) to be only four. 
In \cite{vGW} it
is shown that a desingularization  
$\tilde V_t$ of the generic fiber of the family \eqref{Bertfamily},
obtained by blowing up along its  
singular locus has $h^3 = \dim~H^3(\tilde V_t,\Q) = 6$ ($h^{3,0} = 1 =
h^{0,3}$) and $h^2 = \dim H^2(\tilde V_t, \Q) = h^{1,1} = 138$.

Within the family \eqref{Bertfamily} there are two special fibers with a different 
number of singular points and different Hodge
numbers. The first, is the fiber at 
\[
t = -5 \pm 2\sqrt 5, \quad \ie a = \frac{1}{2},\quad b = 9 \pm
4\sqrt 5.
\]
This
threefold has 126 nodes and it is isomorphic to the one studied by Hirzebruch (\ie the self product 
of a regular pentagon configuration). Its desingularization has $h^3 = 2$ ($h^{2,1} = 0
= h^{1,2}$) and $h^2 = h^{1,1} = 152$ (\cf~\cite{H} and \cite{vGW1}).
It was shown in \cite{vGW} that this variety is 
equivalent to the fiber at $c = \frac{1}{2}$ in the family of threefolds obtained as self-product of the 
family
\begin{equation}\label{regpentagon}
F_c(x,y) = (x+c)(y^4-y^2(2x^2-2x+1)+\frac{1}{5}(x^2+x-1)^2).
\end{equation}

The generic fiber of \eqref{regpentagon} is a skew pentagon but this family satisfies the condition \eqref{condition} only
for $c \in \{\frac{1}{2},-2\}$.
The second special fiber  in \eqref{Bertfamily} is 
obtained for
\[
t = 5 \pm 2\sqrt 5,\quad \ie a = \frac{1\pm \sqrt 5}{2},\quad b =
9\pm 4\sqrt 5
\]
(each couple must be chosen with the same sign). This 
threefold is isomorphic to the fiber at  $c = -2$ in
\eqref{regpentagon}. These varieties have 120 nodes: \ie $10 \cdot 10 + 4 \cdot 4 + 2 \cdot 2$. The two further singular points (with
respect to the singular set of the generic fiber) come up
because $H_t(x,t)$
takes the same value also at the 2 critical points on the x-axis not on the lines (note that 2 among the 10 critical points of $F_{-2}(x,y)$ defined by the intersection points of the lines are on the $x$-axis \cf~\eqref{factorization} below). The
desingularized varieties have $h^3 = 4$ 
($h^{3,0} = h^{0,3} = h^{2,1} = h^{1,2} = 1$) and $h^2 = h^{1,1} =
141$. This paper will focus  on the geometry and
arithmetic of them.

\section{The Chebyshev family of quintic threefold.}\label{1}


In this section we introduce a family of quintic threefolds in
$\Pr^4$ whose geometry (and arithmetic) is related to the 
special fiber at $t = 5 \pm 2\sqrt
5$ in \eqref{Bertfamily}.

We start by introducing the following function in three variables
\begin{equation}\label{pol}
f(y_1,y_2,y_3) := y_1^5 + y_2^5 + y_3^5 + (y_1y_2)^5 + (y_1y_3)^5 + (y_2y_3)^5.
\end{equation}

This polynomial is the  sum of two
Dickson polynomials of the first kind: $D_5^{(i)}(x_1,x_2,1)$, $i = 1,2$. 
We recall, for completeness, the definition of Dickson polynomials
of the first kind in several variables and few of their properties. For a complete report on
this theory we refer to \cite{L}. 

Let $R$ be a commutative ring with identity.

\begin{defn} The Dickson polynomials of the first kind
$D_n^{(i)}(x_1,x_2,\ldots,x_k,a)$, $1 \le i \le k$, $a \in R$ $k \ge
1$ integer, are given by the
functional equation
\[
D_n^{(i)}(x_1,x_2,\ldots,x_k,a) = s_i(y_1^n,\ldots,y_{k+1}^n),\quad 1
\le i \le k,
\]
where $x_i = s_i(y_1,\ldots,y_{k+1})$ and $y_1\cdots y_{k+1} = a$.

The vector $D(k,n,a) = (D_n^{(1)},\ldots,D_n^{(k)})$ of the $k$ Dickson
polynomials is called the Dickson polynomial vector.
\end{defn}

When $k = 2$, the Dickson polynomials can be defined recursively in
the following way
\begin{lem} The polynomials $D_n^{(i)} := D_n^{(i)}(x_1,x_2,a)$, $i =
1,2$ satisfy the recurrence relation
\begin{gather*}
D_n^{(1)} = x_1D_{n-1}^{(1)}- x_2D_{n-2}^{(1)} + aD_{n-3}^{(1)}\\
\text{with}~D_0^{(1)} = 3,~D_1^{(1)} = x_1,~D_2^{(1)} = x_1^2 - 2x_2
\end{gather*}
and
\begin{gather*}
D_n^{(2)} = x_2D_{n-1}^{(2)}- ax_1D_{n-2}^{(2)} + a^2D_{n-3}^{(2)}\\
\text{with}~D_0^{(2)} = 3,~D_1^{(2)} = x_2,~D_2^{(2)} = x_2^2 - 2ax_1.
\end{gather*}
\end{lem}
\begin{proof}  \cf~\cite{L}.\end{proof}

When $k = 2$ and $a = 1$, one has
\begin{gather}\label{Dick}
D_5^{(i)}(x_1,x_2,1) = s_i(y_1^5,y_2^5,y_3^5),\quad 1 \le i \le 2\\
x_i = s_i(y_1,y_2,y_3), \quad\text{and}\quad y_1y_2y_3 = 1\notag\\
s_1 = y_1 + y_2 + y_3,\quad s_2 = y_1y_2 + y_1y_3 + y_2y_3,\quad s_3 =
y_1y_2y_3. \notag
\end{gather}

The sum $D_5^{(1)}(x_1,x_2,1) + D_5^{(2)}(x_1,x_2,1)$ is then 
the polynomial $f$  in \eqref{pol}. 
The recurring relation shows that
$f$ may be described by a symmetric equation of degree $5$ in the
variables $x_1,x_2$ ($x_i$ are the first two elementary symmetric
polynomials in three variables $x_i = s_i(y_1,y_2,y_3)$, under the condition
$y_1y_2y_3 = 1$). It turns out that $f$ is implicitely defined by
\begin{gather*}\label{Cheb1}
f(y_1,y_2,y_3) = D_5^{(1)}(x_1(y_1,y_2,y_3),x_2(y_1,y_2,y_3),1) +
 D_5^{(2)}(x_1(y_1,\ldots),x_2(y_1,\ldots),1) \\ 
D_5^{(1)}(x_1,x_2,1) + D_5^{(2)}(x_1,x_2,1) =
 P_5(x_1,x_2),
\end{gather*}
where $P_5(x_1,x_2)$ is the following symmetric quintic polynomial
\begin{equation}\label{Cheb2}
P_5(x_1,x_2) = (x_1^5+x_2^5) - 5x_1x_2(x_1^2+x_2^2) + 5x_1x_2(x_1+x_2)
+ 5(x_1^2+x_2^2) - 5(x_1+x_2). 
\end{equation} 

Note that \eqref{Cheb2} can be written
as a sum of two determinants (each determinant corresponds to a
Dickson polynomial)
\begin{displaymath}
P_5(x_1,x_2) = \text{det}\left (\begin{array}{ccccc}
x_1 & 1 & 0 & 0 & 0\\
2x_2 & x_1 & 1 & 0 & 0 \\
3 & x_2 & x_1 & 1 & 0 \\
0 & 1 & x_2 & x_1 & 1 \\ 
0 & 0 & 1 & x_2 & x_1
\end{array}\right ) + \text{det}\left (\begin{array}{ccccc}
x_2 & 1 & 0 & 0 & 0\\
2x_1 & x_2 & 1 & 0 & 0 \\
3 & x_1 & x_2 & 1 & 0 \\
0 & 1 & x_1 & x_2 & 1 \\ 
0 & 0 & 1 & x_1 & x_2
\end{array}\right ).
\end{displaymath}

Each determinant represents a generalized Chebyshev polynomial
(of the first kind) associated to the root system $A_2$ (\cf~\cite{HW}). 
The general definition of a Chebyshev polynomial of
the first kind in $k$ indeterminates is given by (\cf~\cite{K})
\begin{gather*}
P_{m,n}^{(-1/2)}(x_1,\ldots,x_k) = \sum_{i=1}^{k+1}\sum_{j=1\atop j\neq
i}u_i^mu_j^{-n}, \quad\text{for}~m,n \in \Z\\
x_i = s_i(y_1,\ldots,y_{k+1}),\quad\text{and}~x_{k+1} = y_1\cdots
y_{k+1} = 1.
\end{gather*}
Hence, one may write $f$  as the
sum $P_{0,-5}^{(-1/2)}(x_1,x_2) +
P_{0,5}^{(-1/2)}(x_1,x_2)$, where the variables $x_1,x_2,x_3$ are
defined as in \eqref{Dick}.
 
Dickson's polynomials share a number of interesting properties,
especially when
the ring $R$ is the finite field ${\mathbb F}_q$ with $q$ elements. We
recall one of these for future reference.

The Dickson polynomial vector $D(k,n,a)$ induces a map of ${\mathbb F}_q^k$ into
itself defined by
\begin{gather*}
D(k,n,a): \mathbb F_q^k \to \mathbb F_q^k \\
(\zeta_1,\ldots,\zeta_k) \mapsto
(D_n^{(1)}(\zeta_1,\ldots,\zeta_k,a),\ldots,D_n^{(k)}(\zeta_1,\ldots,\zeta_k,a)).
\end{gather*}

The following result shows under what condition $D(k,n,0)$
permutes ${\mathbb F}_q^k$.
\begin{prop}\label{permu} $D(k,n,0)$ permutes ${\mathbb F}_q^k$ if and only if
\[
gcd(n,q^s-1) = 1,\quad\text{for}~s = 1,\ldots k.
\]
\end{prop}
\begin{proof}  We refer to \cite{L} Theorem~3.41 for a complete
proof. In \op, it is shown
under which condition the Dickson vector $D(k,n,a)$ induces
permutation on the field ${\mathbb F}_q$, for $a \in {\mathbb F}_q^{\times}$.
\end{proof}

In the rest of this paragraph we will consider the polynomial $P_5(x_1,x_2)$ as a polynomial of degree 5
in the variables $x_1,x_2$, ``forgetting'' 
the description of 
these variables in terms of elementary symmetric polynomials. 
The interpretation of $P_5$ in terms of Dickson's polynomials will be
reconsidered later on in the paper.\vspace{.1in}
  
The function $P_5(x_1,x_2)$ has 16 non
degenerate critical points and its highest homogeneus part is
non-degenerate (\ie it has only one critical point at zero).
A very important property of this polynomial is that of taking only three
different values at its critical points.  One can easily verify that $P_5$ has 10 
critical points with critical value -2 (these occur at the intersection points of the lines), 2 critical points with critical value 6 (they occur on the $x_1$-axis) and 4 critical points with value -3. 
This particular behavior characterizes the
Chebishev polynomials $P_d(x_1,x_2)$ of degree $d$,
independently of the choice of $d$. In \cite{C}, Chmutov 
defines the number of critical points associated to a given critical
value, for any positive d. In this paper we are mainly
interested in the case $d = 5$, so we will skip the general description for
self-contained reasons. \vspace{.1in}

Starting with the Chebyshev polynomial $P_5(x_1,x_2) \in \Q[x_1,x_2]$,
let first define the fibration of affine plane curves (generically
of genus $6$) 
\begin{equation}\label{P5}
\A^2 \to \A^1,\quad P_5(x_1,x_2) = t.
\end{equation}
 
Then, we introduce a new couple of variables $\{x_4,x_5\}$ and
we consider the self-fiber-product of the family \eqref{P5}. We obtain
 the following family of quintic hypersurfaces in $\A^4_{\Q}$
\begin{gather}\label{Cheb}
\A^4 \to \A^1,\quad f_5(x_1,x_2,x_4,x_5) = t\\
f_5(x_1,x_2,x_4,x_5) := P_5(x_1,x_2) - P_5(x_4,x_5).\notag 
\end{gather}

Let $\bar X \subset \Pr^4$ 
be the projective closure of the fiber at $t = 0$
in \eqref{Cheb}. Then, it follows from the definition of (and the value at) the critical
points that
\begin{lem} The projective treefold $\bar X$ has
\[
(10)^2 + (4)^2 + (2)^2 = 120
\]
non-degenerate double points.
\end{lem}\vspace{.1in}
 
Note that the fibers at $t = 5
\pm 2\sqrt 5$ in \eqref{Bertfamily} and  at $c = -2$ in 
 \eqref{regpentagon} are threefolds with the same number of double points. This
analogy is not casual. To explain it, we describe 
more in detail the Chebyshev family. \vspace{.1in}

In the cyclotomic field $\Q(\zeta_5)$ ($\zeta_5$ a $5$-th
primitive root of unit), the contraction-translation of $P_5(x_1,x_2)$:
$-\frac{1}{10}P_5(x_1,x_2)+2$ factorizes as
\begin{multline}
-\frac{1}{10}P_5(x_1,x_2) + 2 =
 -\frac{1}{10}(x_2+\zeta_5+\zeta_5^3-x_1-x_1\zeta_5-x_1\zeta_5^2-x_1\zeta_5^3)(x_2-1-\zeta_5-\zeta_5^3+x_1\zeta_5)\\ (x_2-1-\zeta_5-\zeta_5^2+x_1\zeta_5^2)(x_2+\zeta_5+\zeta_5^2+x_1\zeta_5^3)(x_2+2+x_1).
\end{multline}
In other words, $\Q(\zeta_5)$ is the minimal
field extension of $\Q$ over which every fiber of \eqref{P5} splits
into a product of 5 lines.
Furthermore, the maximal real subfield of the cyclotomic field
$\Q(\zeta_{20})$ is the minimal
field extension over which each fiber of the skew-pentagon
configuration $F_{-2}(x_1,x_2)$ is defined. Let 
$\zeta_{re} = \zeta_{20} + \zeta_{20}^{-1}$. Then, in the extension
$\Q(\zeta_{re})$ the five lines are given by
\begin{multline}\label{factorization}
F_{-2}(x_1,x_2) =
\frac{1}{625}(5x_2+x_1\zeta_{re}^3+5\zeta_{re}-2\zeta_{re}^3)(5x_2+5\zeta_{re}-\zeta_{re}^3-10x_1\zeta_{re}+3x_1\zeta_{re}^3)\\
(5x_2-5\zeta_{re}+\zeta_{re}^3+10x_1\zeta_{re}-3x_1\zeta_{re}^3)(-x_1\zeta_{re}^3-5\zeta_{re}+2\zeta_{re}^3+5x_2)(x_1-2)
\end{multline}

The following proposition shows the existence of an isomorphism
over $\Q(i)$ which identifies the family \eqref{P5} and the skew-pentagon
family \eqref{regpentagon}: $F_{-2}(x_1,x_2) = t$, (for $c = -2$).

\begin{prop}\label{match} Let $F_{-2}(x_1,x_2) =
(x_1-2)(x_2^4-x_2^2(2x_1^2-2x_1+1)+\frac{1}{5}(x_1^2+x_1-1)^2)$ and
$P_5(x_1,x_2) = x_1^5+x_2^5-5(x_1x_2-1)(x_1^2+x_2^2-x_1-x_2)$. Then 
\[
P_5(x_1,x_2) = -10F_{-2}(x_1,x_2)-2\]
by means of the $\Q(i)$-linear map
\[
\A^2_{\Q} \to \A^2_{\Q},\quad (x_1,x_2) \mapsto
(-\frac{x_1+x_2}{2}+1,\frac{i(x_1-x_2)}{2}). 
\] 
\end{prop}
\begin{proof}  It is rather easy to check that the affine 3--folds $F_{-2}(x_1,x_2) - F_{-2}(x_3,x_4) = 0$ and $P_5(x_1,x_2) - P_5(x_3,x_4) = 0$ have the same number of points over fields containing $\sqrt{-1}$. Hence, any bijective linear transformation between $P_5(x_1,x_2)$ and $F_{-2}(x_1,x_2)$ will be defined over such fields. Among the critical points of $F_{-2}(x_1,x_2)$, one finds  the points $(x_1,x_2) = (\frac{2}{3}a^2+\frac{1}{2},a)$ with $a$ satisfying the equation $16a^4-36a^2+9 = 0$. The field generated by these points is $\Q(\sqrt 5,\sqrt 3)$. Similarly, among the critical points of $P_5(x_1,x_2)$ one finds the points $(x_1,x_2) = (-b+\frac{1}{2}-\frac{1}{2}b^3,b)$ with $b$ satisfying the equation $b^4+b^3+2b^2-b+1 = 0$. The field generated by these points is $\Q(\sqrt 5, \sqrt{-3})$. Of course, both the polynomials have further critical points, but we simply look for a linear map that relates the above critical points. This is because the image of points like $(-b+\frac{1}{2}-\frac{1}{2}b^3,b)$ (\cf above) is defined over a subfield of $\Q(\sqrt 5, \sqrt{-3}, i) = \Q(b,i)$ and this field contains $\Q(a)$. A similar argument shows that the image of points like $(\frac{2}{3}a^2+\frac{1}{2},a)$ is defined over a subfield of $\Q(a,i)$ and this field contains $\Q(b)$. One can easily express the parameter $a$ in terms of $b$ and $i$: $a = -\frac{i}{4}b^3 - ib + \frac{i}{4}$. The definition of the linear map follows.\end{proof}\vspace{.1in}

\begin{rem}\label{re1}\end{rem}

Over the field $\Q(\zeta_{20})$ (the common
field over which each line in both constructions is defined),
one may use a cross-ratio method to show that there are only two
possible maps between the 2 sets of lines that define \resp the Chebyshev and the
skew-pentagon families. First, note that 5-tuples of distinct points in general position (\ie no three of
them lying on a line) in $\Pr^2_{\R}$, modulo ${\rm PGL}_3(\R)$, are
the same as 5-tuples of distinct points in $\Pr^1_{\R}$, modulo
${\rm PGL}_2(\R)$.  In $\Pr^1$ one has 2
cross-ratios. Consider the set
$\{(0,0),(0,1),(1,0),(1,1)\}$ as a ``basis'' in $\Pr^1$ and then add the
fifth point $(a,b)$.  In this
way, $a$ and $b$ will play the role of the two invariants. If the five points in
$\Pr^2$ are given as the vectors in $\A^3$: $v_1 =
(v_{1_1},v_{1_2},v_{1_3}),\ldots, v_5 = (v_{5_1},v_{5_2},v_{5_3})$,
then we write the 2 ratio invariants as
\begin{multline*}
r(v_1\vert v_2,\ldots,v_5) :=
\frac{D(v_1,v_2,v_4)D(v_1,v_3,v_5)}{D(v_1,v_2,v_5)D(v_1,v_3,v_4)},\\ 
r(v_2\vert v_3,v_4,v_5,v_1) := 
\frac{D(v_2,v_3,v_5)D(v_2,v_4,v_1)}{D(v_2,v_3,v_1)D(v_2,v_4,v_5)},
\end{multline*}
where
\begin{displaymath}
D(v_i,v_j,v_k) := \text{det}\left (\begin{array}{ccc}
v_{i_1} & v_{i_2} & v_{i_3} \\
v_{j_1} & v_{j_2} & v_{j_3}  \\
v_{k_1} & v_{k_2} & v_{k_3}  \\
\end{array}\right). 
\end{displaymath}\vspace{.1in}

The two 5-tuples of lines associated to the equations $P_5(x,y) = 0 =
F_{-2}(x,y)$, determine 5-tuples of length-3-vectors: $v_1,\ldots,v_5$ and $w_1,\ldots,w_5$. 
We computed explicitely the ratios $[r(v_1\vert v_2,\ldots,v_5), r(v_2\vert v_3,v_4,v_5,v_1)]$
and  we compared them with $[r(w_{\sigma(1)}\vert
w_{\sigma(2)}, \ldots,w_{\sigma(5)}), r(w_{\sigma(2)}\vert
w_{\sigma(3)},w_{\sigma(4)},w_{\sigma(5)},w_{\sigma(1)})]$ for any 
permutation of five elements $\sigma: \{1,2,3,4,5\} \to \{1,2,3,4,5\}$. It turns out that there are
only two possible matchings, and they  correspond  to the two
isomorphisms over $\Q(i)$ of  proposition~\ref{match}.
\newline\rightline{$\Box$}

\section{Group action and cohomology.}\label{2}


In this paragraph we describe a group action on the Chebyshev
threefold $X$ (and more in general on any fiber of \eqref{Cheb}), using which we 
 deduce its Hodge numbers  and the topological Euler
characteristic. 

For resolutions of hypersurfaces singularities in $\Pr^4_{\C}$
 the computation of $h^3$ depends upon the
number of double points s and the defect {\it def}, by means of
the well known formula: $h^3 = 204 -2s + 2\cdot{\it def}$. The defect computes the
number of independent 
divisors on the threefold through the nodes and 
not homologous to a multiple of a generic hyperplane section. Its knowledge depends on 
the special position of the nodes on the embedded threefold.  The defect is 
 a corank of a certain matrix whose
size depends on the degree of the hypersurface and on the
number of its singular points.\vspace{.1in}

We denote by $S$  the 
closed set (for the Zariski topology) described by the singular points
of $X$. We  write $\mathfrak I_S$ for its ideal
sheaf. Let $F_5(x_0,\ldots,x_4)$ be the homogeneus
polynomial whose zeroes define $\bar X$ as a hypersurface in
$\Pr^4$ (we homogenize $f_5(x_1,x_2,x_4,x_5)$ by using a further
variable $x_0$ and then we re-number the variables consecutively).  One
has the following exact sequence of $\C$-vector spaces 
\begin{equation}\label{$*$}
\begin{CD}
0 \to H^0(\Pr^4,\mathfrak I_S\mathcal O(5)) \to H^0(\Pr^4,\mathcal O(5)) @>{\text e}>> H^0(S,\mathcal O_S) 
\to H^1(\Pr^4,\mathfrak I_S\mathcal O(5)) \to 0
\end{CD}
\end{equation}

where
\[
\dim_{\C} H^0(\Pr^4,\mathcal O(5)) =\left( \begin{array}{cc} 9\\5\end{array} \right) = \frac{9!}{5!\cdot 4!} = 126.
\]

The map $e$ evaluates a homogeneus
polynomial $f_i$ of degree five at each point $P_j$ of the
singular locus $S$. It can be described by the  
rectangular matrix $M := (f_i(P_j))$ (whose size is $126 \times 120$) . 
The defect of $\bar X$ is the corank
of $M$ \ie the dimension of $ H^1(\Pr^4,\mathfrak I_S\mathcal O(5))$.

Because  $\dim_{\C}H^0(\Pr^4,\mathfrak I_S\mathcal O(5))
\ge 25$ (any linear combination as $\sum_{i=0}^4
L(x_0,\ldots,x_4)\cdot \frac{\partial F_5}{\partial
x_i}(x_0,\ldots,x_4)$ vanishes at $S$, where $L$ a linear homogeneus
polynomial in five variables), one necessarily expects that 
\[
\text{corank}~M = \dim_{\C}H^1(\Pr^4,\mathfrak I_S\mathcal O(5)) \ge
19.
\]

The evaluation of the rank of $H^0(\Pr^4,\mathfrak I_S\mathcal O(5))$, depends 
upon a precise understanding of  the special position of the nodes of
$X$ in $\A^4_{\C}$. Note that $\bar X - X$
is a non singular variety.\vspace{.1in}  

There are three involutions $\sigma_i$ ($i = 1,\ldots,3$)
acting on the affine threefold $X$ and more in general on each fiber of the family \eqref{Cheb}. They are 
\begin{equation}\label{inv}
\begin{split}
\sigma_1(x_1,x_2,x_3,x_4) =
(x_2,x_1,x_3,x_4),&\quad\sigma_2(x_1,x_2,x_3,x_4) = 
(x_1,x_2,x_4,x_3),\\
\sigma_3(x_1,x_2,x_3,x_4) &= (x_3,x_4,x_1,x_2).
\end{split} 
\end{equation} 
The group generated by them is the dihedral group of order eight
$\mathfrak D_4$,
that defines the rotations and the reflexions of an affine plane
preserving a regular polygon with four vertices. Let $\sigma_5 = \sigma_1\sigma_3$. In terms of the the involutions \eqref{inv}, $\mathfrak D_4$ has the following presentation
\[
\mathfrak D_4 = <\sigma_3, \sigma_5~\vert~\sigma_3^2 = \sigma_5^4 = 1,~\sigma_3\sigma_5\sigma_3 = \sigma_5^{-1}>.
\]

It is easy to deduce from this description that $\mathfrak D_4$ may be realized as a semidirect product $<\sigma_3> \ltimes <\sigma_5>$. The element $\sigma_5$  
generates a cyclic group of order four that describes the rotations of the plane through an angle of $\frac{\pi}{2}$. We indicate this group by $\mathfrak C_4$. The
element $\sigma_3$ (one of the four reflections of the plane) generates a cyclic group of order two and it conjugates elements of $\mathfrak C_4$ into their inverses.
The center $Z(\mathfrak D_4)$ of $\mathfrak D_4$ is generated by $\sigma_5^2$ ($= \sigma_1\sigma_2$).
  
The exact sequence \eqref{$*$} becomes a sequence of
$\mathfrak D_4$--modules. In the rest of this paragraph we provide a
descritpion of the map $e$, restricted to each isotypic space for
the action of the characters of $\mathfrak D_4$.  

It is a standard fact that each element of $\mathfrak D_4$ can be written uniquely either in the form $\sigma_5^k$ or in the form $\sigma_3\sigma_5^k$ for $0\le k \le 3$. We define $\sigma_4=\sigma_5^2$, $\sigma_6 = \sigma_5^3$, $\sigma_7=\sigma_3\sigma_5^2$, whereas we have that $\sigma_1=\sigma_3\sigma_5^3$ and $\sigma_2=\sigma_3\sigma_5$.

The group $\mathfrak D_4$ has five conjugacy classes and therefore five
characters. Each of the four 1-dimensional irreducible characters $\chi_i$, $i =
1,\ldots,4$, factorizes through the center $Z(\mathfrak D_4)$ (\ie
$\chi_i(\sigma_5^2) = 1$). These characters are 
uniquely characterized by their values on $\sigma_1$ and
$\sigma_3$. 
The fifth character $\chi_5$ is obtained by inducing a 1-dimensional
irreducible representation of $\mathfrak C_4$. Let $\chi$ be its character. We choose for it the following description
\[
\chi: \mathfrak C_4 \to \C^*,\quad \chi(\sigma_5^k) = \exp(\frac{\pi ik}{2}),\quad 0 \le k \le 3.
\]

The element $\sigma_3$ acts on $\chi$ as 
$(\sigma_3\chi)(\sigma_5^k) = \chi(\sigma_3\sigma_5^k\sigma_3)$. Because this action  does not stabilize the
character $\chi$ ($\sigma_3\chi \neq \chi$), one deduces that 
the identity  is the only irreducible representation of the 
isotropy group of $\chi$.
The definition of  the
character $\chi_5$ of the corresponding irreducible 2-dimensional
representation of $\mathfrak D_4$ follows then by  a standard induction
argument. Explicitely, we have  
\[
\chi_5 = \text{Ind}_{\mathfrak C_4}^{\mathfrak D_4}(\chi);\quad \chi_5(\sigma_5^k) =
2\cos\frac{\pi k}{2},\quad \chi_5(\sigma_3\sigma_5^k) = 0,\quad 0 \le k \le 3.
\]

The following table resumes the above description.
\begin{eqnarray}\label{ta}
\begin{tabular}{|c|c|c|c|c|c|}
characts/elts & $\chi_1 = id$ & $\chi_2$ & $\chi_3$ & $\chi_4$ & $\chi_5$\\ \hline\hline
id & $1$ & $1$ & $1$ & $1$ & $2$ \\ \hline
$\sigma_3\sigma_5^3$ & $1$ & $-1$ & $1$ & $-1$ & $0$\\ \hline
$\sigma_3\sigma_5$ & $1$ & $-1$ & $1$ & $-1$ & $0$\\ \hline
$\sigma_3$ & $1$ & $1$ & $-1$ & $-1$ & $0$\\ \hline
$\sigma_5^2$ & $1$ & $1$ & $1$ & $1$ & $-2$\\ \hline
$\sigma_5$ & $1$ & $-1$ & $-1$ & $1$ & $0$\\ \hline
$\sigma_5^3$ & $1$ & $-1$ & $-1$ & $1$ & $0$\\ \hline
$\sigma_3\sigma_5^2$ & $1$ & $1$ & $-1$ & $-1$ & $0$\\ \hline
 \end{tabular}
\end{eqnarray}


Next, we study the behavior of the evaluation map $e$
on the isotypic subspaces. Our goal is
the determination of the rank of the vector spaces $ker(e) \simeq H^0(\Pr^4,\mathfrak I_S\mathcal O(5))$ and $coker(e) \simeq H^1(\Pr^4,\mathfrak I_S\mathcal O(5))$, together with  
the description of their $\mathfrak D_4$-module structure. This will be accomplished by a
direct computation of the rank of their isotypic
subspaces.  

As a linear map of vector spaces, $e$ is
described by a rectangular matrix of size $126\times 120$. Because $e$
is $\mathfrak D_4$--linear,  it carries the  
isotypic subspace for a certain irreducible representation to the
subspace relative to the same representation within 
the vector space of the singular points. 

 The space $H^0(\Pr^4,\mathcal O(5))$ is generated by 
monomials of degree five in the homogeneus variables
$x_0,\ldots,x_4$. On each of the five affine charts  
that cover the projective space $\Pr^4$, one writes down the 126
monomials in the corresponding  
four affine coordinates. The group $\mathfrak D_4$ acts on them by switching the
coordinates $x_i$ (\cf~\eqref{inv}).
The matrix that describes
this action has size $126\times 126$. One may consider its
projection onto the  eigenspace  
associated to each irreducible character of $\mathfrak D_4$.
 
Let $\rho: \mathfrak D_4 \to \bf{GL}$$(H^0(\Pr^4,\mathcal O(5))$ denote the linear 
representation of $\mathfrak D_4$, with character $\phi$ and let $V_i$ ($i =
1,\ldots,5$) be the isotypic  
spaces corresponding to the irriducible characters $\chi_i$ of degrees $n_i$. 
The projection formula onto the i-th isotypic space, reads as 
\[
p_i = \frac{n_i}{8}\sum_{\sigma\in \mathfrak D_4}\chi_i(\sigma)\rho(\sigma).
\]

One writes down the canonical decomposition of 
$H^0(\Pr^4,\mathcal O(5))$ into irreducible representations, each of which 
appears with the corresponding multiplicity
\begin{equation}\label{mult}
(\phi\vert \chi_i) = <\phi,\chi_i> = \frac{1}{8}\sum_{\sigma\in
\mathfrak D_4}\phi(\sigma) \chi_i(\sigma) = \frac{1}{8}\sum_{\sigma\in
\mathfrak D_4} \chi_i(\sigma)Tr(\rho(\sigma)).
\end{equation}

The computation of the multiplicities in \eqref{mult} can be easily done
 by hands. The idea is to keep track of the number of monomials that are
 invariant under the action of the group $\mathfrak D_4$. The
 explicit description of the eigenvectors generating 
 each isotypic subspace was obtained with the
 help of a computer. One gets the
following decomposition
\begin{prop} The canonical decomposition of $H^0(\Pr^4,\mathcal
O(5))$ as a $\mathfrak D_4$-module is
\[
H^0(\Pr^4,\mathcal O(5)) = 27V_1 \oplus 9V_2 \oplus 23V_3 \oplus 7V_4 
\oplus 30V_5.
\]
We denoted by $V_i$ the irreducible representation associated to the
character $\chi_i$. The number in front of each eigenspace corresponds
to its degree.  
\end{prop}

To get a similar decomposition for the vector
space of the singular points, one has to 
determine their coordinates explicitely 
and understand the orbits for the action of $\mathfrak D_4$ on them.
We accomplished this part using some invariant properties of the Chebyshev polynomial $P_5$.
We review them briefly in the following.
 
Let think of $P_5(x_1,x_2)$ as a map $P_5: \C^2_{(x_1,x_2)} \to \C$.
We consider a set of coordinates $\{u_1, u_2, u_3\}$ satisfying the
equation $u_1 + u_2 + u_3 = 0$. These coordinates are related to
$x_1$ and $x_2$ by the map $h$
\[
h: \C^2_{(u_1,u_2,u_3)} \to \C^2_{(x_1,x_2)}, 
\]
\begin{equation}\label{h}
\begin{split}
x_1 &= \exp(-2\pi iu_1) + \exp(-2\pi iu_2) + \exp(-2\pi iu_3)\\
x_2 &= \exp(-2\pi
i(u_1+u_2)) + \exp(-2\pi i(u_1+u_3)) + \exp(-2\pi i(u_2+u_3)).
\end{split}
\end{equation}\vspace{.1in}

The function $h$ plays for the generalized Chebyshev polynomial  a
 role analogous to that played by the cosine  
for the classical Chebyshev polynomial in one variable. Namely, $h$
is invariant under the affine Weyl group associated to the root system
$A_2$. We
introduce a second set of coordinates $\{u_4, u_5, u_6\}$ satisfying
relations similar to those described by the first set. The set
$\{u_4, u_5, u_6\}$  
is related to the coordinates $x_4,x_5$ via equations analogous to \eqref{h}.  
Then
\begin{equation}\label{acca}
P_5(h(u_1,u_2,u_3)) - P_5(h(u_4,u_5,u_6)) = 0,\quad u_1+u_2+u_3 = 0 =
u_4 + u_5 + u_6 
\end{equation}

describes an affine covering $Y$ of the threefold $X$.
To make as simple as possible the notations, we rename the complex
variables as 
\[
y_1 = \exp(-2\pi iu_1),~y_2 =  \exp(-2\pi iu_2),~y_4 = \exp(-2\pi
iu_4),~y_5 =  \exp(-2\pi iu_5).
\]  

Then,  $Y$ is defined in terms of the variables $y_i$, 
as the set of solutions of
\begin{equation}\label{sols}
\begin{split}
y_1^5 + y_2^5 + y_3^5 + y_1^{-5} + y_2^{-5} + y_3^{-5} &= y_4^5 +
y_5^5 + y_6^5 + y_4^{-5} + y_5^{-5} + y_6^{-5}  \\
y_3 = (y_1y_2)^{-1}, &\quad y_6 = (y_4y_5)^{-1}.
\end{split}
\end{equation}

Namely, $Y$ is defined in terms of  the polynomials $f$ introduced in 
\eqref{pol} as the sum of two Dickson
polynomials of the first kind. 

The algebraic description of the covering $\pi: Y \to X$ 
is induced by the restriction of the map 
\begin{gather}\label{map}
\pi^*: \C[x_1,x_2,x_3,x_4,x_5,x_6] \to \C[y_1,y_2,y_3,y_4,y_5,y_6]\\
\pi^*(x_1,x_2,x_3,x_4,x_5,x_6) = 
(y_1+y_2+y_3,y_1y_2+y_1y_3+y_2y_3,y_1y_2y_3,\notag\\
y_4+y_5+y_6,y_4y_5+y_4y_6+y_5y_6,y_4y_5y_6)\notag  
\end{gather}

to the four dimensional affine space $\A^4 = \{x_3 = 1 = x_6\}$. One obtains
\[
\pi^*: \C[x_1,x_2,x_4,x_5]/(f_5(x_1,x_2,x_4,x_5)) \to
 \C[y_1,y_2,\ldots,y_6]/(y_1y_2y_3-1,y_4y_5y_6-1)
\]

where $f_5$ was defined in \eqref{Cheb}. If $P =
(y_1,\ldots,y_6)$ is a point in $Y$, the definition of
the map $\pi: Y \to X$ is given by
\begin{align*}
\pi(P) &= \pi(y_1,\ldots,y_6) =\\
&= (y_1+y_2+\frac{1}{y_1y_2},y_1y_2+\frac{1}{y_1}+\frac{1}{y_2},
y_4+y_5+\frac{1}{y_4y_5},y_4y_5+\frac{1}{y_4}+\frac{1}{y_5}).
\end{align*}

The quotient $X$ is the variety of the orbits for the action of
$(\mathfrak S_3)^2$ on $Y$ ($\mathfrak S_i =$ symmetric
group in $i$ variables). Each of the two factors $\mathfrak S_3$ acts on the
sets of coordinates $\{y_1,y_2,y_3\}$ \resp $\{y_4,y_5,y_6\}$. The fiber of
$\pi$ over a 
point $Q = (x_1,x_2,x_4,x_5) \in X$, coincides with its (closed) orbit
and it is given by the set (\ie closed zero scheme) of the ordered pairs
$(y_1,y_2,y_3,y_4,y_5,y_6)$ of points of $Y$ (at most 36 distinct,
multiple pairs are possible),
where $y_1,y_2,y_3$ \resp $y_4,y_5,y_6$ range among the
permutations of the zeroes of the equation
\begin{equation}\label{eq1}
T^3 -x_1T^2+x_2T-1 = 0;\quad\text{for:} 
\quad \begin{array}{c} x_1 =
y_1+y_2+y_3 \\
x_2 = y_1y_2+y_1y_3+y_2y_3 \\
y_1y_2y_3 = 1 
\end{array} 
\end{equation}
\resp
\begin{equation}\label{eq2}
T^3 -x_4T^2+x_5T-1 = 0;\quad\text{for:} 
\quad \begin{array}{c} x_4 =
y_4+y_5+y_6 \\
x_5 = y_4y_5+y_4y_6+y_5y_6 \\
y_4y_5y_6 = 1. 
\end{array} 
\end{equation}

The coordinate ring of $X$
coincides with the $\mathfrak S_3^2$-invariants of the coordinate ring
of $Y$: \ie the homomorphism \eqref{map} is described by the embedding $A(X)
= A(Y)^{\mathfrak S_3^2} \hookrightarrow A(Y)$.
This agrees with the fact that $f_5$ is a linear combination of
elementary symmetric functions in the two separate sets of variables
$\{x_1,x_2\}$ and $\{x_4,x_5\}$.\vspace{.1in} 



Using the equality $\cos(x) =
\frac{e^{ix} + e^{-ix}}{2}$ one gets
\[
P_5(h(u_1,u_2,u_3)) = 2\cos(10\pi u_1) + 2\cos(10\pi u_2) + 2\cos(10\pi
u_3).
\]

By a direct computation (\cf~\cite{C}), it is easy to see that if 
$(\zeta_1,\zeta_2,\zeta_3)$ is a critical
point of the compositum $P_5(h(u_1,u_2,u_3))$, then $(\zeta_i)^{30} = 1$, $i =
1,\ldots,3$. The values taken at the critical points are
$6,~-2,~-3$. In other words, the $6$-tuple
$(\zeta_1,\zeta_2,\zeta_3,\zeta_4,\zeta_5,\zeta_6)$ is a singular
point of $Y$ if (and only if) it satisfies the conditions  
\begin{equation}\label{cond}
\begin{split}
\zeta_i = &e^{-\pi i\frac{\alpha_i}{15}},\quad \alpha_i \in \Z,~ i =
1,2,4,5 \\
&\zeta_3 = (\zeta_1\zeta_2)^{-1},\quad \zeta_6 = (\zeta_4\zeta_5)^{-1}.
\end{split}
\end{equation}

\begin{prop} The singular locus of the quasi-projective variety $Y$
consists of $8750$ double points.
\end{prop}
\begin{proof}  If $P$ is a point in the singular locus of the variety
$Y$ (\cf~\eqref{sols}), then the coordinates  of $P$ must satisfy the
conditions \eqref{cond}. This implies that for $P =
(\zeta_1,\ldots,\zeta_6) \in Sing(Y)$, $\zeta_i^5 = e^{-\pi
i\frac{\alpha_i}{3}}$ for $i = 1,2,4,5$ and $\alpha_i =
1,\ldots,5$. Therefore, the singular locus of $Y$ is
described by the set of points $(e^{-\pi
i\frac{\alpha_1}{3}}, e^{-\pi i\frac{\alpha_2}{3}}, e^{-\pi
i\frac{(\alpha_1+\alpha_2)}{3}}, e^{-\pi
i\frac{\alpha_4}{3}}, e^{-\pi i\frac{\alpha_5}{3}}, e^{-\pi i\frac{(\alpha_4+\alpha_5)}{3}})$ in the singular
locus of the variety
\begin{equation}\label{ram}
x_1 + x_2 + x_1x_2 +  x_1^{-1} + x_2^{-1} + (x_1x_2)^{-1} = x_4 +
x_5 + x_4x_5 +  x_4^{-1} + x_5^{-1} +  (x_4x_5)^{-1}.
\end{equation}
The $4$-tuples $(\alpha_1,\alpha_2,\alpha_4,\alpha_5)$ are
chosen among the set of solutions of
\[
\cos(\pi\frac{\alpha_1}{3}) + \cos(\pi\frac{\alpha_2}{3}) +
\cos(\pi\frac{\alpha_1+\alpha_2}{3}) = \cos(\pi\frac{\alpha_4}{3}) +
\cos(\pi\frac{\alpha_5}{3}) + \cos(\pi\frac{\alpha_4+\alpha_5}{3})
\]

under the condition $0 \le \alpha_i \le 5$.
By a direct computation one obtains the following locus
\begin{multline*}
T = \{(1,1,1,1,1,1), (1,-1,-1,1,-1,-1), (1,-1,-1,-1,1,-1),
(1,-1,-1,-1,\\-1,1), (e^{-2\pi i/3},
 e^{-2\pi i/3},e^{-2\pi i/3},e^{-2\pi i/3},e^{-2\pi i/3},e^{-2\pi
i/3}), (e^{-2\pi i/3},e^{-2\pi i/3},e^{-2\pi i/3},e^{2\pi
i/3},\\e^{2\pi i/3},e^{2\pi i/3}), (-1,1,-1,1,-1,-1),
(-1,1,-1,-1,1,-1), (-1,1,-1,-1,-1,1),\\(-1,-1,1,1,-1,-1),
(-1,-1,1,-1,1,-1), (-1,-1,1,-1,-1,1), \\(e^{2\pi i/3},
 e^{2\pi i/3},e^{2\pi i/3},e^{-2\pi i/3},e^{-2\pi i/3},e^{-2\pi
i/3}), (e^{2\pi i/3},
 e^{2\pi i/3},e^{2\pi i/3},e^{2\pi i/3},e^{2\pi i/3},e^{2\pi
i/3})\}.
\end{multline*}

Hence, to each of the 14 points in $T$
correspond $5^4$ points solutions of \eqref{sols}. In this way, one
counts globally $5^4 \times 14 = 8750$ singular points for $Y$. One
can direct verify that they are all rational double points.
\end{proof}\vspace{.1in}

The semi-direct product $G := (\mathfrak S_3^2 \rtimes \mathfrak S_2) \ltimes
\mu_5^4$ acts on $Y$ ($\mu_5 =$ fifth roots of
unity). The action of $\mathfrak S_2$ interchanges the
two sets of coordinates $\{y_1,y_2,y_3\}$ and $\{y_4,y_5,y_6\}$ and it
is preserved by the morphism $\pi$. It corresponds to the action of
the involution $\sigma_3$ on $X$ (\cf\eqref{inv}). The
abelian subgroup $\mu_5^4$ acts freely on 
$Y$ and on its singular set. The action of  $\mathfrak S_3
\times \mathfrak S_3$ has instead stabilizers. If $S'
= Sing(Y)$, then the action of $\mu_5^4$ on 
$S'$ determines the 14 
orbits of points described by the set $T$.
In other words, the vector space $H^0(S',\mathcal O_{S'})$, inherits a
structure of $\mu_5^4$--module. The related $\mu_5^4$-representation
is a direct sum of 14 copies of the regular representation
$H^0(S',\mathcal O_{S'}) = \bigoplus_{14} \C[\mu_5^4]$.

We divide the 8750 singular points in $S'$ in accord with
the action that $\mathfrak S_3 \times \mathfrak S_3$ has on them. We are mainly
interested in understanding the unramified part of this action 
(\ie to focus on the points not contained in the ramification locus
\eqref{ram} of $\pi$). In fact, the unramified points define a
set (\ie zero scheme) on which $\mathfrak S_3^2$ acts freely: the
quotient space for this action determines the (non singular)
zero-scheme $S$ of $X$-rational double points. Each
point in $T$ has a non-trivial stabilizer for the action of $\mathcal
S_3 \times \mathfrak S_3$ and therefore the image by $\pi$ of points
within their $\mu_5^4$-orbits determines the singular locus of $X$.

In the following part we give some details on how we accomplished these calculations: the
numbers were checked with the help of a computer.
Within the following five orbits of $T$
\begin{multline*}
\{(1,1,1,1,1,1), (e^{-2\pi i/3},
 e^{-2\pi i/3}, e^{-2\pi i/3}, e^{-2\pi i/3}, e^{-2\pi i/3}, e^{-2\pi
i/3}), (e^{-2\pi i/3}, e^{-2\pi i/3},\\ e^{-2\pi i/3}, e^{2\pi
i/3}, e^{2\pi i/3}, e^{2\pi i/3}), (e^{2\pi i/3},
 e^{2\pi i/3}, e^{2\pi i/3}, e^{-2\pi i/3}, e^{-2\pi i/3}, e^{-2\pi
i/3}), (e^{2\pi i/3}, e^{2\pi i/3},\\ e^{2\pi i/3}, e^{2\pi i/3},
 e^{2\pi i/3}, e^{2\pi i/3})\}
\end{multline*}
one counts $5\cdot (12)^2 = 720$
$(\mathfrak S_3 \times \mathfrak S_3)$-unramified $Y$-points. More
precisely, each of these orbits
contains representatives of $144/36 = 4$ $(\mathfrak S_3 \times
\mathfrak S_3)$-(unramified) orbits. The complete orbits
 determine $4 \times 5 = 20$ $X$-singular points globally.  In the
remaining $9$ orbits of $T$  
there are $9\cdot (20)^2 = 3600$ $(\mathfrak S_3 \times
\mathfrak S_3)$-unramified $Y$-points subdivided in the following way.  
Four among the $9$ orbits contain $8$
representatives of $(\mathfrak S_3 \times 
\mathfrak S_3)$-(unramified) orbits each. The (complete) orbits associated
to these representatives (that are distributed in different
$\mu_5^4$-orbits of the set $T$)  determine
$32$ $X$-singular points globally. 
Four further orbits (among the 9) contain 16 representatives of
$(\mathfrak S_3 \times \mathfrak S_3)$-(unramified) orbits each. The
associated complete orbits 
determine 64 singular points on $X$. 
Finally, the last $\mu_5^4$-orbit contains 4 representatives of
$(\mathfrak S_3 \times \mathfrak S_3)$-(unramified) orbits. Hence one counts
$4$ singular points on $X$. 
The total number of $(\mathfrak S_3 \times \mathfrak S_3)$-unramified
$Y$-points is $720 + 3600 = 4320$. They correspond  to 
$4320/36 = 120$ singular points of $S$.  

Using the explicit description of these points in logarithmic
$(u)$-coordinates (\cf~\eqref{h}), one determines the $\mathfrak D_4$-action
on them. First, one writes down the
character of the representation $\tilde\rho: \mathfrak D_4 \to
\bf{GL}$$(H^0(S,\mathcal O_S))$ as a trace of the automorphism given
by the action of $\mathfrak D_4$ on the 120 $(\mathfrak S_3 \times \mathfrak
S_3)$-unramified points in the $14$ orbits of $T$. Then, one computes
the multiplicities as in \eqref{mult} of each
irreducible representation associated to $\tilde\rho$. 
This can be easily done with the help of a computer. We obtained the following
canonical decomposition
\begin{prop} The canonical decomposition of the $\mathfrak D_4$-module
$H^0(S,\mathcal O_S)$ is
\[
H^0(S,\mathcal O_S) = 27U_1 \oplus 13U_2 \oplus 17U_3 \oplus 7U_4
\oplus 28U_5.
\]
Here $U_i$ is the irreducible representation associated to the
character $\chi_i$ and the number in front of it corresponds to its
multiplicity.
\end{prop}\vspace{.1in}

For simplicity, we re-number the projective coordinates used in
\eqref{map} as $\{x_0,x_1,x_2,x_3,x_4\}$. 

Next, we describe the decomposition in eigenspaces of the
$25$-dimensional subvector space $K \subset ker(e)$, generated by the linear
combinations 
\[
\sum_{i=0}^4L(x_0,\ldots,x_4)\cdot
\frac{\partial F_5}{\partial x_i}(x_0,\ldots,x_4).
\]

This is obtained using techniques similar to the ones already described.
Firstly, one writes the character associated
to the representation $\mathfrak D_4 \to GL(K)$ and then one computes the
multiplicities of each irreducible representation associated to it.
The action of $\mathfrak D_4$ on the five partial derivatives 
and on the linear functions $L$
is given by
\[
\sigma_i(\frac{\partial F_5}{\partial x_j}) = \sum_j
a_{ij}\sigma_i(\frac{\partial F_5}{\partial x_j}),
\quad\text{and}\quad 
\sigma_i(x_j) = \sum_j b_{ij}\sigma_i(x_j).
\]

The character associated to this representation
is the tensor product of the two
characters $\chi(\sigma_i) = \text{tr}(a_{ij}(\sigma_i))$ and
$\chi'(\sigma_i) = \text{tr}(b_{ij}(\sigma_i))$. 
Straightforward computations determine the following table
\begin{eqnarray*}
\begin{tabular}{|c|c|c|}
characts/elts & $\chi$ & $\chi'$ \\ \hline\hline
id & $5$ & $5$ \\ \hline
$\sigma_1$ & $3$ & $3$\\ \hline
$\sigma_2$ & $3$ & $3$\\ \hline
$\sigma_3$ & $-1$ & $1$\\ \hline
$\sigma_4$ & $1$ & $1$\\ \hline
$\sigma_5$ & $-1$ & $1$\\ \hline
$\sigma_6$ & $-1$ & $1$\\ \hline
$\sigma_7$ & $-1$ & $1$\\ \hline
 \end{tabular}
\end{eqnarray*}

The values taken by the tensor product $\chi\otimes \chi'$
(\ie the character of the 25--dimensional representation) are readable
from the above table.

The scalar products
\[
m_i = <\chi\otimes\chi',\chi_i> =
\frac{1}{8}\sum_{\sigma_j\in
\mathfrak D_4}(\chi\otimes\chi')(\sigma_j)\cdot\chi_i(\sigma_j)
\]

determine the multiplicites of each irreducible representation. One obtains
\[
K = 5K_1 \oplus K_2 \oplus 6K_3 \oplus K_4 \oplus 6K_5.
\]

It remains to describe the canonical decomposition and the dimension of
$ker(e)$ itself.
This can be done by evaluating the eigenvectors that define
the  eigenspaces associated to each irreducible character at the 120
double points. For each of 
these five maps we computed its kernel
and  rank and  obtained  the following canonical decomposition
\begin{prop} $H^0(\Pr^4,\mathfrak I_S\mathcal O(5)) \simeq Ker(e) = 5K_1 \oplus K_2 \oplus (6K_3 \oplus \C v_3) \oplus K_4
\oplus 6K_5$. 
\end{prop}

The following corollary resumes what we obtained.
\begin{cor} With the above notations, one has
\[
rk~H^0(\Pr^4,\mathfrak I_S\mathcal O(5)) = 26,\quad rk~
H^1(\Pr^4,\mathfrak I_S\mathcal O(5)) = 20.
\]
In particular, the canonical decomposition of the $\mathfrak D_4$--module $H^1(\Pr^4,\mathfrak I_S\mathcal O(5))$ is
\[
H^1(\Pr^4,\mathfrak I_S\mathcal O(5)) = 5U_1' \oplus 5U_2' + U_3' + U_4' + 4U_5'
\]
where $U_i'$ is the irreducible representation associated to the character $\chi_i$.
\end{cor} 

\begin{rem}\label{re2}\end{rem}

For the description of the
eigenvector $v_3 \in ker(e)- K$ 
one looks at the the complement of $6K_3$ in the
23--dimensional 
isotypic space $V_3 \subset H^0(\Pr^4,\mathcal O(5))$ associated to
the character $\chi_3$. One searches for a linear
combination of the $17$ elements generating the complement of 
$6K_3$ in the space $23V_3$. The
following $\chi_3$--invariant affine polynomial satisfies these
conditions and it was found with the help of a computer.
\begin{multline*}
h(x_1,\ldots,x_4) =
-x_1^2x_3-x_2^2x_3+x_1x_3^2+x_2x_3^2-x_1^2x_4-x_2^2x_4+x_1x_4^2+x_2x_4^2+(-x_1^4x_3-\\-x_2^4x_3+x_1x_3^4+x_2x_3^4-x_1^4x_4-x_2^4x_4+x_1x_4^4+x_2x_4^4)
+
[-x_1x_2x_3-x_1x_2x_4+x_1x_3x_4+\\+x_2x_3x_4+(-x_1^2x_2^2x_3-x_1^2x_2^2x_4+x_1x_3^2x_4^2+x_2x_3^2x_4^2)]+[-x_1^3x_3-x_2^3x_3+x_1x_3^3+x_2x_3^3-x_1^3x_4-\\-x_2^3x_4+x_1x_4^3+x_2x_4^3-(x_1^2x_2x_3-x_1x_2^2x_3-x_1^2x_2x_4-x_1x_2^2x_4+x_1x_3^2x_4+x_2x_3^2x_4+\\+x_1x_3x_4^2+x_2x_3x_4^2)]
\end{multline*} \newline\rightline{$\Box$}\vspace{.1in}

From these computations one may deduce the
Hodge structure on $H^*(\tilde X,\C)$. We are only concerned with the knowledge of  $H^2(\tilde X)$ and
$H^3(\tilde X)$ as one knows that $h^1(\tilde X) = 0$. The following
complex of sheaves in $\Pr^4$ 
\begin{equation}\label{com}
\begin{CD}
\mathcal E^\cdot:\quad 0 @>>> 0 @>>> 0 @>>> \Omega^3_{\Pr^4}(X) @>d>> \mathfrak 
I_S\omega_{\Pr^4}(2X)
\end{CD}
\end{equation}  

induces in hypercohomology the isomorphism
\begin{equation*}
F^2H^i(\tilde X) \simeq H^{i+1}(\Pr^4,\mathcal E^\cdot).
\end{equation*}
For a proof of this statement we refer to \cite{S}. The only
interesting result is related to the index $i = 3$. Namely one gets
\begin{equation}\label{iso}
F^2H^3(\tilde X,\C) \simeq \frac{H^0(\Pr^4,\mathfrak
I_S\omega_{\Pr^4}(2X))}{\text{Image}(H^0(\Pr^4,\Omega^3(X)) \to
H^0(\Pr^4,\mathfrak I_S\omega_{\Pr^4}(2X)))}.
\end{equation}

Hence

\begin{cor}\label{cor} Under the same notations
\[
h^{3,0}(\tilde X) = 1 = h^{2,1}(\tilde X).
\]
\end{cor}
\begin{proof}  It follows from the above isomorphism and from the fact
that $h^0(\Pr^4,\Omega^3(X)) = 24$ (independent of the number of the
singular points). For smooth, quintic hypersurfaces in $\Pr^4$ one knows that
$h^{3,0} = 1$.
\end{proof}
\begin{rem}\label{re3}\end{rem}

The $\chi_3$--invariant polynomial $h$ defined in Remark~\ref{re2}
is a generator of the space 
$H^{2,1}(\tilde X)$. From \eqref{com} and \eqref{iso} one
obtains  the diagram
\begin{equation*}
\begin{CD}
0 @>>> H^0(\Pr^4,\Omega^3(X)) @>d>> H^0(\Pr^4,\mathfrak
I_S\omega_{\Pr^4}(2X)) @>\beta>> F^2H^3(\tilde X,\C) @>>> 0 \\
@. @. @A{\simeq}A{\psi}A @. @. \\
 @. @. H^0(\Pr^4,\mathfrak I_S\mathcal O(5)) 
\end{CD}
\end{equation*}

The isomorphism $\psi$ is induced by the isomorphism of sheaves
$\omega_{\Pr^4} \simeq \mathcal O_{\Pr^4}(-5)$. The map $\beta$ is
deduced from the blow-up of $X$ along $S$ and by the Poincar\'e residue. The
composite of these maps is described by 
\[
P(x_0,\ldots,x_4) \mapsto \frac{P(x_0,\ldots,x_4)\cdot
\sum_{i=0}^4(-1)^idx_0\wedge\ldots\wedge\hat{dx_i}\wedge\ldots\wedge
dx_4}{F_5^2(x_0,\ldots,x_4)}.
\]

Here $P(x_0,\ldots,x_4) \in H^0(\Pr^4,\mathfrak I_S\mathcal O(5))$ is a
homogeneus polynomial of degree $5$ that vanishes at each point of $S$
and $F_5(x_0,\ldots,x_4)$ is the homogeneus polynomial that
defines $X$.

It is worth to remark that if one chooses $P(x_0,\ldots,x_4) =
 F_5(x_0,\ldots,x_4)$, then its image by means of $\beta\cdot\psi$
gives a generator of the subvector space $H^0(\tilde
 X,\omega_{\tilde X})$ of the holomorphic $3$-forms  on $\tilde
 X$. In other words, the K\"ahler differential from $\Pr^4$ via
 Poincar\'e residue pulls-back to a regular $3$-form on $\tilde X$. If
 instead one chooses  $P = h$, then it follows from  
Remark~\ref{re2} that the image of it via $\beta\cdot\psi$
generates the $(2,1)$-part of the $H^3(\tilde X)$. 

Finally, let   $y_i = \frac{x_i}{x_0}$
(a similar argument works in each of the four affine charts), then 
the differential form $\omega = dy_1\wedge dy_2\wedge dy_3\wedge
dy_4$ is $\chi_2$--invariant as $\sigma_1(\omega) = \sigma_2(\omega) =
-\omega$ and $\sigma_3(\omega) = \omega$ (\cf the character table
description).  Therefore, the image of any
$\chi_3$--invariant polynomial (\eg $F_5$ and $h$)
 by $\beta\cdot\psi$, becomes $\chi_2\chi_3$ ($=
\chi_4$)-invariant.  Because   
$\chi_4(\sigma_1) = -1 = \chi_4(\sigma_3)$, it follows that one cannot
split the rank $4$ motive associated to $F^2H^3(\tilde X)$ by the
action of any subgroup of $\mathfrak D_4$ on $\tilde X$.\newline\rightline{$\Box$}\vspace{.1in}

We recall the well known fact
\begin{lem} Let $\bar X$ be an hypersurface of degree five in $\Pr^4$ with
$120$ nodes (\ie rational double points). Then, the Euler
characteristic of a desingularization $\tilde X$ of $\bar X$ obtained by
blowing up the nodes is given by the formula
\[
\chi(\tilde X) = \chi(X_{\eta}) + 4\cdot 120 = -200 + 4\cdot 120 = 280.
\] 
Here $\bar X_{\eta}$ is the generic fiber of a Lefschetz fibration of
degree five hypersurfaces in $\Pr^4$ with $\bar X$ as a special fiber.
\end{lem}
\begin{proof}  The vanishing cycle exact sequence applied to a Lefschetz
pencil of quintic hypersurfaces in $\Pr^4$ with degenerating fiber $\bar X$
and generic fiber $\bar X_{\eta}$, reads as
\[
0 \to H^3(\bar X,\Q) \to H^3(\bar X_{\eta},\Q) \to V \to H^4(\bar X,\Q) \to
H^4(\bar X_{\eta},\Q) \to 0
\]
\[
H^i(\bar X,\Q) \simeq H^i(\bar X_{\eta},\Q),\quad\text{for}~i \neq 3,4.
\]
Here, $V$ is the $\Q$-vector space of vanishing cycles. Because $\bar X$
has only double points, the rank of $V$ 
equals the number of double points of $\bar X$. Hence one gets
\begin{equation}\label{chi}
\chi(\bar X) = \chi(\bar X_{\eta}) + 120.
\end{equation}
The Euler charactersitic of $\bar X_{\eta}$ (a non-singular hypersurface of
degree $5$ in $\Pr^4$) is given by integrating its third Chern class
against $\bar X_{\eta}$ and then by interpreting the integral using Stokes theorem
\[
\chi(\bar X_{\eta}) = \int_{\bar X_{\eta}}c_3(\bar X_{\eta}) = \int_{\Pr^4}[5J]\wedge
c_3(\mathcal O_{\Pr^4}).
\]
Here, $J = c_1(\mathcal O_{\Pr^4}(1))$. It is well known that for
degree $q$ smooth hypersurfaces in $\Pr^n_{\C}$ the r-th Chern class
is given by
\[
c_r = [\sum_{k=0}^r\binom{n+1}{k}(-q)^{r-k}]J^r.
\]
For degree $5$ hypersurfaces in $\Pr^4$
one gets $\chi(\bar X_{\eta}) = -200$. Using Lefschetz hyperplane theorem, one may deduce that $h^2(\bar X_{\eta}) =
h^{1,1}(\bar X_{\eta}) = 1$ and $h^1(\bar X_{\eta}) = 0$. Hence,
$h^{2,1}(\bar X_{\eta}) = 101$ as 
$h^{3,0}(\bar X_{\eta}) = 1$. In turn, this result implies $\chi(\bar X)
= 2 + h^2(\bar X) - h^3(\bar X) + h^4(\bar X) = 2 + 1 - h^3(\bar X) +
h^4(\bar X)  = 3 - h^3(\bar X) + h^4(\bar X) = -80$.  Hence
$h^3(\bar X) - h^4(\bar X) = 83$. The Leray spectral sequence
for blow-up of double points $\tilde X \to \bar X$  gives the following
exact sequences 
\begin{equation}\label{a2}
\begin{split}
0 \to H^2(\bar X,&\Q) \to H^2(\tilde X,\Q) \to H^2(E,\Q) \to H^3(\bar X,\Q) \to
H^3(\tilde X,\Q) \to 0\\
&0 \to H^4(\bar X,\Q) \to H^4(\tilde X,\Q) \to H^4(E,\Q) \to 0
\end{split}
\end{equation}
\[
H^i(\bar X,\Q) = H^i(\tilde X,\Q),\quad\text{for}~i \neq 2,3,4.
\]
Here $E$ is the exceptional fiber, \ie a union of $120$ $\Pr^1
\times \Pr^1$. From this sequence it follows that $\chi(\tilde X) =
\chi(\bar X) + 120\cdot(\chi(\Pr^1\times\Pr^1) - 1)$. Using \eqref{chi} one concludes that
\[
\chi(\tilde X) = \chi(\bar X_{\eta}) + 120\cdot\chi(\Pr^1\times\Pr^1) = -200 + 4\cdot 120 = 280.
\]
\end{proof}
As an immediate consequence we have
\begin{cor}\label{dimensions} $h^{1,1}(\tilde X) = 141 = h^2(\tilde X),\quad h^3(\bar X) =
104,~\text{and}\quad h^4(\bar X) = 21$.
\end{cor}
\begin{proof}  The equality $h^2(\tilde X) = h^{1,1}(\tilde X)$ follows from
 \eqref{a2} and the fact that the rank-one vector space $H^2(\bar X,\Q)$
is generated by the class of a hyperplane section. The equality
$h^2(\tilde X) = 141$ is a direct consequence of the above lemma and
corollary~3.7 (\ie $h^3(\tilde X) = 4$). The equality $h^3(\bar X)
 = 104$ follows from \eqref{a2}. Finally, $h^4(\bar X) =
21$ is deduced from the difference  $h^3(\bar X) - h^4(\bar X) = 83$
(\cf~proof of the last lemma).  
\end{proof}

\section{Reduction mod. $p$ and counting points.}\label{J1}


In this paragraph we describe the reduction of $\tilde X$ modulo
a prime. We use the same notation as was 
 introduced in the previous paragraphs.

\begin{lem}\label{tablelemma}
The threefold $\tilde X$ has good reduction outside 2, 3 and 5.
\end{lem}
\begin{proof}
In order to prove this, we study the behaviour of the critical points of
$P_5(x_1,x_2)$ mod $p$. In paragraph~\ref{1}, we observed that
$P_5(x_1,x_2)$ has $10$ critical points with value $-2$, $4$ with value $-3$
and $2$ with value $6$. We list these points in the following table.
\begin{eqnarray}
\label{crittable}
\begin{tabular}{|c|c|c|c|c|c|}
\hline
point & \# Gal.orbit & field of def. & value \\
\hline
\hline
$\left(w,w\right)$ & 2 & $\Q(w)$ & $6$ \\
$\left(w+1,w+1\right)$ & 2 & $\Q(w)$ & $-2$ \\
$\left(-\zeta_5,-\zeta_5^{-1}\right)$ & 4 & $\Q(\zeta_5)$ & $-2$ \\
$\left(\zeta_5-\zeta_5^{-1}-1,\zeta_5^{-1}-\zeta_5-1\right)$ & 4 
& $\Q(\zeta_5)$ & $-2$ \\
$\left(w\zeta_3,w\zeta_3^{-1}\right)$ & 4 & $\Q(w,\zeta_3)$ & $-3$ \\
\hline
\end{tabular}
\end{eqnarray}
\vspace{.1in}
One can easily check that these $16$ critical points reduce to different
points \mod primes of $\Q(\zeta_{15})$ above rational primes bigger
than $5$. One way to verify this is to compute the discriminant with respect to
$x_2$ of the resultant of $\frac{dP_5}{dx_1}$ and $\frac{dP_5}{dx_2}$
with respect to $x_1$. This number is only divisible by the primes:
$2$, $3$, $5$, $11$, $19$ and $31$, so these are the only primes mod which 
the critical points can coincide. A finite calculation shows that the
critical points reduce to different points mod primes above $11$, $19$
and $31$.

For $p\neq 5$,
the resultant of $\frac{dP_5}{dx_1}$ and 
$\frac{dP_5}{dx_2}$ (say with respect to $x_1$) is not identically
zero. This implies that these polynomials
do not have a common factor modulo $p$. By Bezout's
theorem one concludes that $P_5(x_1,x_2)$ \mod~$p$ has at most $16$ critical points.
Hence, for $p>5$ there are precisely $16$ critical points:
 the reductions of the $16$ critical points in
characteristic zero.

Note that the values $-2$, $-3$ and $6$ that $P_5$ assumes at the
critical points remain different when reduced modulo primes bigger
than $5$. It follows that the number of singular points of the
affine singular threefold defined by $P_5(x_1,x_2) = P_5(x_4,x_5)$
does not increase upon reducing \mod primes bigger than $5$.

At each of the critical points of $P_5$ one can check that 
the local expansion of $P_5$ has a non-degenerate degree-$2$-part
which remains non-degenerate when one reduces it \mod primes of 
$\Q(\zeta_{15})$, above rational primes bigger than $5$. It follows
that the part of $\tilde X$ above the affine part 
$P_5(x_1,x_2) = P_5(x_4,x_5)$ remains non-singular after the reduction
\mod~$p>5$.

It remains to verify that $\tilde X$ \mod~$p$ has no singularities 
at infinity. In fact, the homogenenized $P_5(x_1,x_2) - P_5(x_4,x_5)$ has the 
following form: $x_1^5+x_2^5-x_4^5-x_5^5+z\,q(x_1,x_2,x_4,x_5,z)$. 
If its partial derivatives were all zero at a point 
where $z=0$ then, in characteristic $\neq 5$, one would have 
$x_1 = x_2 = x_4 = x_5 = 0$. \end{proof}

We are interested in counting $\F_q$-rational points on $\tilde X$.
This will be accomplished in three steps:

\begin{itemize}
\item Determine the number of points on the affine singular part
\item Determine how many points are added in the blow-up
\item Determine the number of points at infinity.
\end{itemize}

An efficient way for computing the number of $\F_q$-rational points 
on the affinesingular part
is as follows. Let $x_1$ and $x_2$ run through $\F_q$, and count
how many times $P_5$ assumes each value. Clearly, the number of
points  will be 
the sum of the squares of these values.

In the blow-up each singular point gets replaced by its 
projectivized tangent cone. Geometrically this is a 
$\Pr^1\times\Pr^1$, but a priori, there could be different Galois
actions on it. In our case it turns out that the rulings of
the cone are defined over the same field as the field of definition
of the singular point. Therefore, each $\F_q$-rational singular point
contributes $(q+1)^2$ points in the blow-up (instead than  $1$
point on the singular model). Using the list of critical points
of $P_5$ as in the above table, one finds that the number of $\F_q$-rational
singular points is $120$ if $q\equiv 1$ \mod~$15$, $104$ if $q\equiv 11$
\mod~$15$, $24$ if $q\equiv 4$ \mod~$15$, $8$ if $q\equiv 14$ \mod~$15$ and
zero otherwise.

The locus at infinity is defined by the homogeneous equation
$x_1^5 + x_2^5 - x_4^5 - x_5^5$. One can determine the number of 
solutions by adding the squares of the number times 
$x_1^5+x_2^5$ assumes each value. In this way one gets the number of solutions
of the homogeneous equation. One has to substract $1$, and divide
by $q-1$ in order to obtain the number of points.

It turns out that for many $q$ one does not need a computer
for counting points. One can use the following proposition 
instead.

\begin{prop}
\label{count}
If $q\equiv \pm 2\ \mod\ 5$ and $\gcd(q,30)=1$,
 then $\tilde X(\F_q)=q^3+q^2+q+1$.
\end{prop}

\begin{proof}
Suppose $q\equiv 2{\rm\ or\ }3\ \mod\ 5$. 
In paragraph~\ref{0} we showed that 
the polynomial $P_5(x_1,x_2) = D_5^{(1)} + D_5^{(2)}$, with
$D_5^{(1)}$ and $D_5^{(2)}$ Dickson
polynomials. It follows from proposition~\ref{permu} 
that the map $(D_5^{(1)},D_5^{(2)})$ permutes $\F_q\times\F_q$.
Hence $P_5\ \mod\ p$ assumes each value exactly
$q$ times and therefore there are $\sum_1^qq^2=q^3$ solutions to the
equation $P_5(x_1,x_2) = P_5(x_4,x_5)$ in $F_q$.
From table~\ref{crittable} it follows that there are no critical points
defined over $\F_q$, hence there is no extra contribution coming
from the blow-up.

The map $x\mapsto x^5$ permutes the elements of $\F_q$, hence
$x_1^5+x_2^5$ assumes each value $q$ times. So 
$x_1^5+x_2^5 = x_4^5+x_5^5$ has $q^3$ solutions in $F_q$, 
which correspond to $\frac{q^3-1}{q-1}$ points.
In total, one counts globally $q^3+\frac{q^3-1}{q-1}=q^3+q^2+q+1$
points.
\end{proof}

\begin{rem}
The determination of the number of points
at infinity also holds for $q\equiv 4\ \mod\ 5$.
\end{rem}

\section{Computing the $L$-function.}\label{J2}


We are interested in computing the $L$-function of the 
Galois action on $H^3(\tilde X_{\overline\Q},\Ql)$.
We start by recalling its definition.
Let $I_p$ be an inertia group at $p$, and $\fr_p$ a
(geometric)
Frobenius element of $\Gal_\Q$. 
The $L$-function is defined as the product
$$L(s)=\prod_{p\rm\ prime}L_p(s)$$
with $$L_p(s)=\left.\det(\fr_p\cdot T-\id)^{-1}\right|_{T=p^{-s}}.$$
Here $(\fr_p\cdot T-\id)$ acts on the $I_p$-invariant subspace of
$H^3(\tilde X_{\overline\Q},\Ql)$.

At the primes $p$ where $\tilde X$ has good reduction, $L_p(s)$
can be determined by counting the points on $\tilde X$ over various
finite fields, and using the Lefschetz trace formula. We
will explain this fact here in more detail. We refer to section~\ref{J4} for the discussion on how to obtain
likely candidates for the
$L_p(s)$ at the primes $p$ where $\tilde X$ has bad reduction.

Suppose that $\tilde X$ has good reduction at a prime $p$. Then, for
$q=p^n$, the Lefschetz trace formula is
$$\tilde X(\F_q)=\sum_{i=0}^6 (-1)^i\tr(\fr_q|H^i(\tilde X_{
\overline\Q},\Ql)).$$
One knows that the trace of the Frobenius on $H^0$ is 1, and on $H^6$ is $q^3$.
The $H^1$ and $H^5$ vanish. The $H^2$ and $H^4$ are algebraic, and
they are related by Poincar\'e duality.
They have dimension 141 (\cf corollary~\ref{dimensions}),
so there exists an integer $k$, with $-141<k<141$ such that
$\tr(\fr_q|H^2(\tilde X_{\overline\Q},\Ql))=k\,q$ and
$\tr(\fr_q|H^4(\tilde X_{\overline\Q},\Ql))=k\,q^2$.
Hence with the the trace formula one obtains
\begin{eqnarray}
\tr(\fr_q|H^3(\tilde X_{\overline\Q},\Ql))=
1+q^3+k(q+q^2)-\tilde X(\F_q).\label{tracefor}
\end{eqnarray}

Therefore, the computation of the trace on $H^3$ is accomplished once one knows $k$ and
$\tilde X(\F_q)$. In the previous section it was explained
how to compute $\tilde X(\F_q)$. To determine $k$ one could try
to understand the Galois action on the divisors generating 
 $H^2$. However, we preferred to proceed differently. In paragraph~\ref{2} we proved that $\dim H^3=4$ (\cf corollary~\ref{cor}),
and by the Weil conjectures one knows that the eigenvalues of 
$\fr_q|H^3$ have absolute value $q^{3/2}$. Hence one gets the inequality
\begin{eqnarray}   
\left|1+q^3+k(q^2+q)-\tilde X(\F_q)\right|
\le 4q^{3/2}.\label{weilest}
\end{eqnarray}
For $q$ big enough, there will be a unique $k$ that satisfies
this inequality. We found that $k$ is uniquely determined if $q>20$.

The characteristic polynomial of the $\fr_q$ action on $H^3$
can be expressed in terms of the traces of $\fr_q$ and 
$\fr_{q^2}$. Let $a_q=\tr(\fr_q|H^3)$, then the characteristic
polynomial is
$$T^4-a_qT^3+\frac{1}{2}(a_q^2-a_{q^2})T^2-q^3a_qT+q^6,$$
and for $q=p$ the local $L$-factor becomes
$$L_p(s)=\frac{1}
{1-a_pp^{-s}+\frac{1}{2}(a_p^2-a_{p^2})p^{-2s}-a_pp^{3-3s}+p^{6-4s}}.$$
This explains how to compute $L_p$ for
$p>20$.

For $7\le p\le 19$, we determined $a_{p^2}$, $a_{p^3}$ and
$a_{p^4}$. The first and third number enables one to 
compute the eigenvalues of $\fr_{p^2}$. The
eigenvalues of $\fr_p$ are the square roots of these numbers.
 To determine
which square root one should take, one can use the value of $a_{p^3}$.
The first few traces  are listed in the following table
\begin{eqnarray}
\begin{tabular}{|c||c|c|c|c|c|c|c|}
\hline
$p$          & 7 &   11 & 13   & 17   &  19  & 23      \\
\hline
$a_p$        & 0 &$-116$&  0   & 0    &$-20$ & 0       \\
\hline
$a_{p^2}$ &$-140$& 1444 & 5980 &$-340$& 6404 & 6900  \\
\hline
\hline
$p$       &  29    &   31   & 37     & 41     & 47     & 53  \\
\hline
$a_p$     &  60    &   24   & 0      & $-316$ &  0     &  0  \\
\hline
$a_{p^2}$ &$-95116$&$-82876$&$-59940$&$-51516$& 187060 &$-471700$ \\
\hline
\hline
$p$        &      59  &      61  &   67      &  71     &  73  &  79  \\
\hline
$a_p$      & $-1160$  & $-1116$  &    0      &  $-156$ &  0   &$-460$ \\
\hline
$a_{p^2}$  &$-146156$ &$-131436$ &$-907180$ & $-814316$&27740 &$-1520396$\\
\hline

\end{tabular}\label{tracetable}
\end{eqnarray}

\section{Two 2-dimensional representations.}\label{J3}


It is important to remark that the traces $a_q$ listed in the previous 
table are $0$ for
$q\equiv\pm 2\ \mod\ 5$. In this paragraph we will prove this and 
we will show
that  this behaviour implies that the Galois representation is 
induced.

\begin{lem}
\label{trace0}
If $q\equiv\pm 2\ \mod\ 5$ then the trace $a_q$ is 0.
\end{lem}

\begin{proof}
It follows from propostion $\ref{count}$ that
$\#\tilde X(\F_q)=q^3+q^2+q+1$. The inequality in $\eqref{weilest}$ yields:
$$\left|(k-1)(q^2+q)\right|\leq4q^{3/2}.$$
Hence, $k=1$ for $q$ big enough. Using the trace formula
$\eqref{tracefor}$ it follows that $a_q=0$. For small $q$ the 
$a_q$ are listed in the table $\eqref{tracetable}$.
\end{proof}

Let $\rho$ denote the semisimplification of the Galois representation
$$\Gal({\overline\Q}/\Q)\longrightarrow
\Aut\left(H^3(\tilde X_{\overline\Q},\Ql)\otimes\Ql(\sqrt{5})\right).$$
First we will discuss
some properties of the characteristic polynomials of the 
restriction
$\rho|_{\Gal({\overline\Q}/F)}$. If $p$ is inert in the extension
$F/\Q$, and $\wp$ is the prime above $p$ then $\fr_\wp=(\fr_p)^2$,
and the characteristic polynomial of $\fr_\wp$ is
$$(T^2-\frac{1}{2}a_{p^2}T+p^6)^2.$$
If $p$ is split, and $\wp$ is a prime above $p$, then $\fr_\wp$ is
a Frobenius element at $p$, and hence $\fr_\wp$ and $\fr_p$ have the
same characteristic polynomial.

Clearly, the characteristic polynomials of 
$\rho|_{\Gal({\overline\Q}/F)}$ at inert primes 
split in two quadratic factors. The polynomials we
computed at split primes also split in two quadratic factors
over $\Ql(\sqrt{5})$. For example, at $p=11$ the polynomial is
$$\left(T^2+(58+2\sqrt{5})T+1331\right)
\left(T^2+(58-2\sqrt{5})T+1331\right).  $$
It turns out that there is no
quartic field such that all the characteristic polynomials split
further in this quartic field.
This suggests the following theorem.

\begin{thm}\label{Jaspth}
The restriction $\rho|_{\Gal({\overline\Q}/F)}$ is reducible. It is a 
direct sum of two 2-dimensional representations $\sigma$ and $\sigma'$.
One has $\rho=\Ind_F^\Q\sigma$.
\end{thm}

\begin{proof}
Let $V$ denote the vectorspace 
$H^3(\tilde X_{\overline\Q},\Ql)\otimes{\overline\Ql}$, and let
$\tilde\rho:\Gal({\overline\Q}/\Q)\rightarrow\Aut(V)$ be the
extension of scalars of $\rho$. 
Let $\chi$ be the Dirichlet character of $\Gal({\overline\Q}/\Q)$
with kernel $\Gal({\overline\Q}/F)$. By lemma $\ref{trace0}$ 
we have that $\tr(\tilde\rho(\fr_p))=0$ for all $\fr_p$ not in 
$\ker\chi$.
So $\tr(\tilde\rho(\fr_p))=\tr((\tilde\rho\otimes\chi)(\fr_p))$ 
for all $p>5$. It
follows that $\tilde\rho$ and $\tilde\rho\otimes\chi$ are isomorphic. Let
$T:V\rightarrow V$ be an isomorphism that intertwines $\tilde\rho$ and
$\tilde\rho\otimes\chi$. So $T(\tilde\rho(g)v)=\chi(g)\tilde\rho(g)T(v)$
for all $g\in\Gal({\overline\Q}/\Q)$ and $v\in V$. 
Let $\lambda$ be an eigenvalue of $T$, and $V_\lambda$ its eigenspace.
Clearly, $T$ can not be a scalar map, so $V_\lambda$ is a proper
subspace of $V$. Since $\tilde\rho$ and $\tilde\rho\otimes\chi$ 
actually have coefficients in $\Ql$, the eigenvalue $\lambda$ is 
contained in an extension of $\Ql$ of degree at most 4.
The image $\tilde\rho(\Gal({\overline\Q}/F))$ commutes with $T$,
so its action on $V_\lambda$ is a subrepresentation of
$\tilde\rho|_{\Gal({\overline\Q}/F)}$. From the remarks made above
about the splitting behaviour of the characteristic polynomials it
follows that the only possible dimension of $V_\lambda$ is 2, and
that $\lambda$ is contained in $\Ql(\sqrt{5})$. 
It also follows 
$T$ has only one other eigenspace $V_{\lambda'}$, with eigenvalue 
$\lambda'$.
Denote the representation on $V_\lambda$ by $\sigma$ and the 
representation on $V_{\lambda'}$ by $\sigma'$.

Let $g\in\Gal({\overline\Q}/\Q)$ be an element that represents the
non-identity element of $\Gal(F/\Q)$.
Suppose
$v\in V_\lambda$. Then
$$T(\tilde\rho(g)v)=-\tilde\rho(g)(T(v))=-\lambda\tilde\rho(g)v.$$
So $\tilde\rho(g)v$ is an eigenvector with eigenvalue 
$-\lambda=\lambda'$, and
$\tilde\rho(g)$ sends $V_\lambda$ to $V_{\lambda'}$. It follows
that $\rho=\Ind_F^\Q\sigma$.
\end{proof}

We end this paragraph with a discussion on what could be a
different reason for the fact that 
$\rho$ is induced. 

Remember that the variety $\tilde X$
was obtained as desingularisation of the closure of an affine
variety $X$. And $X$ is a quotient of a variety $Y$ defined
by the equations
\begin{eqnarray*} 
\begin{tabular}{c}
$y_1^5+y_2^5+y_3^5+y_1^{-5}+y_2^{-5}+y_3^{-5}=
y_4^5+y_5^5+y_6^5+y_4^{-5}+y_5^{-5}+y_6^{-5}$,\\
$y_1y_2y_3=1,\qquad y_4y_5y_6=1$.
\end{tabular}
\end{eqnarray*}
We quotient by the action of the group $\mathfrak S_3\times \mathfrak S_3$, where the
first $\mathfrak S_3$ permutes the coordinates $y_1$, $y_2$ and $y_3$, and the
second permutes $y_4$, $y_5$ and $y_6$.

On $Y$ there are a lot of automorphisms defined over $\Q(\zeta_5)$.
For example, one can send 
$(y_1,y_2,y_3)$ to $(\zeta_5y_1,\zeta_5^{-1}y_2,y_3)$, and do something
similar with $(y_4,y_5,y_6)$. 
One can push down these automorphisms to obtain
correspondences on $X$. It can easily be shown that pushing down an
automorphism or its $\Q(\zeta_5)/\Q(\zeta_5+\zeta_5^{-1})$-conjugate
results in the same correspondence on $X$. So the correspondences one
gets this way are defined over $\Q(\zeta_5+\zeta_5^{-1})=F$.
They give rise to linear endomorphisms on the
cohomology that commute with the $\Gal({\overline\Q}/F)$-action.
So the eigenspaces of these endomorphisms are subrepresentations
of $\tilde\rho|_{\Gal({\overline\Q/F})}$. 
We believe that they are 2-dimensional, and that they are infact
the representations $\sigma$ and $\sigma'$.

\section{The $L$-function.}\label{J4}


Conjecturally, the $L$-function has an analytic continuation
to $\C$, and satisfies a functional equation (cf. \cite{Se}).
In this section we explain how one may use this to find candidates
for the local $L$-factors at the bad primes.

Let $N$ be the conductor of $\rho$, as defined in \cite{Se}.
It is not known whether $N$ is independent of $\ell$.
The complete $L$-function (cf. \cite{Se}) has the following shape
$$\Lambda(s)=N^{s/2}(s\pi)^{-2s}\Gamma(s)\Gamma(s-1)L(s).$$
The conjectured functional equation of $\Lambda$ is
\begin{equation}\label{functeq}
\Lambda(s)=w\Lambda(4-s)
\end{equation}
for some $w\in\{\pm1\}$.

The product $\prod L_p(s)$ is known to converge for 
$\Re(s)\ge\frac{5}{2}$. In what it follows we use a well known trick that allows one to test 
numerically the functional equation under the assumption of its
analytic continuation (see e.g. \cite{BST}).
Let $m\geq0$ be an integer such that
all the derivatives $L^{(j)}(2)$ vanish for $0\leq j\leq m-1$.
Choose a real number $t$. 
The idea is to integrate the function of one complex variable 
$$G(s):= \frac{\Lambda(2+s)t^{-s}}{2\pi i\,s^{m+1}}$$
along a path in $\C$ around 0.
By Cauchy's theorem
this integral is $\Lambda^{(m)}(2)/(m!)$.
By stretching the path more and more in such a way that in the
limit it becomes a union of two lines, (the first going from $r-i\infty$ to 
$r+i\infty$ and the other going from $-r+i\infty$ to $-r-i\infty$,
for some real $r>\frac{1}{2}$), one obtains
\begin{eqnarray}
\frac{\Lambda^{(m)}(2)}{m!}
=\frac{1}{2\pi i}\int_{r-i\infty}^{r+i\infty}
\frac{\Lambda(2+s)t^{-s}}{s^{m+1}}\,ds+\frac{1}{2\pi i}
\int_{-r+i\infty}^{-r-i\infty}\frac{\Lambda(2+s)t^{-s}}{s^{m+1}}ds.
\label{int1}
\end{eqnarray}

Define
$$F_m(x)=\frac{1}{2\pi\,i}\int_{r-i\infty}^{r+i\infty}
\frac{\Gamma(s+2)\Gamma(s)}{x^{s}s^m}ds.$$
We refer to \cite{BST} for a discussion on how to compute
values of $F_m$ efficiently.

The first integral of $(\ref{int1})$ is equal to
\begin{eqnarray*}
&\displaystyle
\frac{1}{2\pi i}\int_{r-i\infty}^{r+i\infty}
\Gamma(s+2)\frac{\Gamma(s+1)}{s^{m+1}}
\left(\frac{4\pi^2}{\sqrt{N}}\right)^{-s-2}t^{-s}
\sum_{n=1}^{\infty}\frac{a_n}{n^{s+2}}ds\ =\\
&\displaystyle
\frac{N}{16\pi^4}\sum_{n=1}^{\infty}\frac{a_n}{n^2}\frac{1}{2\pi i}
\int_{r-i\infty}^{r+i\infty}\Gamma(s+2)\Gamma(s)\,
\left(\frac{4nt\pi^2}{\sqrt{N}}\right)^{-s}ds\ =\\
&\displaystyle
\frac{N}{16\pi^4}\sum_{n=1}^{\infty}\frac{a_n}{n^2}
F\left(\frac{4nt\pi^2}{\sqrt{N}}\right).
\end{eqnarray*}

Using the functional equation
\eqref{functeq} one can rewrite the second integral of 
\eqref{int1} as
$$\frac{w}{2\pi i}\int_{r-i\infty}^{r+i\infty}
\frac{\Lambda(2+s)t^s}{s^{m+1}}
\ =\ \frac{wN}{16\pi^4}\sum_{n=1}^{\infty}\frac{a_n}{n^2}
F_m\left(\frac{4n\pi^2}{t\sqrt{N}}\right).$$

Consequently, one obtains
\begin{eqnarray}
L^{(m)}(2)\ =\ m!\sum_{n=1}^{\infty}\frac{a_n}{n^2}
F_m\left(\frac{4nt\pi^2}{\sqrt{N}}\right)
\ +\ wm!\,\sum_{n=1}^{\infty}\frac{a_n}{n^2}
F_m\left(\frac{4n\pi^2}{t\sqrt{N}}\right).\label{12}
\end{eqnarray}

We know how to compute $L_p(s)$ for $p>5$. But we don't know how to
compute $N$, $w$, and $L_p(s)$ for $p\leq 5$. The idea is that we make
a guess for these. 
For each guess one can use
\eqref{12} to compute $L^{(m)}(2)$ for a few $m$. We did this
computation for several values of $t$. One knows that $L^{(m)}(2)$
is independent of $t$. Hence, if the computations give different 
results for different $t$, one has probably made wrong
guesses. We made new guesses until the computation appeared to
be independent of $t$.
The candidates for $w$ are $\pm1$. Because the 
threefold $\tilde X$ has bad reduction
only at $2$, $3$ and $5$, these are the only primes that
can appear in the factorization of $N$. Hence, 
one may assume that $N = 2^a3^b5^c$ for
some integers $a$, $b$ and $c$.
To guess $L_p(s)$ for $p\leq 5$
we recall that the characteristic polynomial of Frobenius acting on 
the inertia invariants has degree at most 4.
Its roots are conjectured to have absolute value $p^{j/2}$ for some $0\leq j\leq 3$
(cf. \cite{Se}).
This limits the possible choices to a finite number. One can reduce
this number by using the fact that $\rho$ is
induced from a 2-dimensional representation $\sigma$. It is
well-known that inducing a Galois representation does not change
the $L$-function. Hence, one should find the characteristic polynomials
of $\sigma$. These polynomials have degree at most $2$, and if the representation
really ramifies at some $p$ (i.e. $p|N$), their  degrees are at most
$1$. In this case there are only $9$ possibilities to be tested, 
namely the polynomials $1$ and
$T\pm (N\wp)^{j}$ for $0\leq j\leq 3$. The corresponding local
$L$-factors are $1$ and $(1\pm p^{2j-2s})^{-1}$ for $p=2$ or $3$, and
$1$ and $(1\pm 5^{j-s})^{-1}$ for $p=5$. 

For $m=0$, $w=-1$, $N=2^23^25^4$ and
$L_2(s)=(1-2^{2-2s})^{-1}$, $L_3(s)=(1-3^{2-2s})^{-1}$ and
$L_5(s)=1$ we computed $(\ref{int1})$ by considering all terms with 
$n<30000$ in the infinite sums. We repeated this computation for $t=1$, $1.4$,
$1.6$, $2.0$ and $2.5$. We found the same value  
$0$ upto $48$ decimal places. This clearly suggests that the $L$-series
vanishes at $s=2$.
For $m=1$ we computed $(\ref{int1})$ for $n<3000$. We took the same
values for $t$, and found that the answer remained constant upto
14 decimal places. 
The value for $L'(s)$ we found in this way is 
$$L'(s)=2.83811389801282$$

\begin{rem}
The exponents in the conductor we found are equal to the codimension
of the inertia invariants implied by the degree of the local
$L$-factors. This suggests that $\rho$ is only tamely ramified away
from $\ell$.
\end{rem}

\section{The description of a quaternion algebra and an Eichler order.}\label{3}


This paragraph is devoted to the definition of a quaternion algebra and an
order. These objects will be used later  in the description of the
Brandt matrices. The computation of the class and type number of the order are also included.
The following notations will hold throughout the rest of the paper.

We denote by $F$ the totally real field $\Q(\sqrt 5)$ and by  
$w = \frac{1+\sqrt 5}{2}$ a fundamental unit. We write $\mathfrak
o_F$ for the ring of integers of $F$. For $p \in \Z$ a rational prime
ideal, we write $\wp_p$ for a prime ideal of $\mathfrak o_F$ above $p$.

Let $B$ be the unique (up to $F$--linear isomorphism) 
definite quaternion algebra
with center $F$,
ramified at both the infinite places  and at $\wp_2$ and $\wp_3$ (\ie with
discriminant $D(B/F) = \wp_2\wp_3$).

\begin{lem}\label{quat} The quaternion algebra $B$ is uniquely
determined by the following setting
\[
B =  L + LY,\quad L = F(X)
\]
where
\[
X^2 = -6,\quad Y^2 =  \sqrt 5~\frac{1-\sqrt 5}{2} 
\]
for $X,Y \in B$ and $Y\cdot l = \bar l\cdot Y,~\forall l \in L$.
Here, $l \mapsto \bar l$ denotes the non trivial $F$-automorphism of
 $L$.
\end{lem}

In characteristic different from two, $B$ can be
shorthly described as
\[
B = (-6,\sqrt 5~\frac{1-\sqrt 5}{2}) = (-6,w-3).
\]
$B$ is the central, simple, $4$--dimensional division algebra over
$F$ generated by the elements $1,X,Y,XY \in B$ satisfying the
following relations
\[
X^2 = -6,\qquad Y^2 = \sqrt 5~\frac{1-\sqrt 5}{2},\qquad XY = -YX.
\]

\begin{proof}  We only show that $B$  ramifies at $\wp_2, \wp_3$ and at both
the infinite places. For more details on the properties of $B$, as
well as for the proof of its uniqueness up to isomorphism, we refer to
\cite{V}. Notice that $B$ is totally definite (\ie $B \otimes_{\Q} \R$ is
a division quaternion algebra isomorphic to two copies of the Hamilton's
quaternions). In fact, $-6$ and $\sqrt 5~\frac{1-\sqrt 5}{2}$ are
totally negative numbers. It remains to show that $\wp_2,\wp_3$ are the two
finite primes where $B$ ramifies. In the local field $F_{\wp_3}$,
the Hilbert symbol (= Hasse invariant) is given by 
\[
(-6,\sqrt 5~\frac{1-\sqrt 5}{2})_{\wp_3} = (\frac{\sqrt
5~\frac{1-\sqrt 5}{2}}{3})_{Lg} = -1,
\]
where  the subscript $Lg$ denotes the Legendre
symbol. In fact,  $\sqrt 5~\frac{1-\sqrt 5}{2}$ is not a square \mod $\wp_3$.
Hence $B_{\wp_3} = B \otimes_F F_{\wp_3}$ is a field. In other words, at the inert prime ideal $\wp_3 \in F$, the
algebra $B$ ramifies. At the ramified 
prime $\wp_5$ (ramified for the extension $F/\Q$), the Hilbert symbol is
\[
(-6,\sqrt 5~\frac{1-\sqrt 5}{2})_{\wp_5} = 1.
\]
Hence, the place $\wp_5$ does not ramify in $B$. 
It follows from the
definition of the algebra that the only 
possible place where $B$ could possibly ramify is $\wp_2$. This is in
fact the case as a quaternion algebra is supposed to ramify at an even
number of primes in a field. 
\end{proof}

It follows from Lemma~\ref{quat} that one can define on $B$ an
involutive anti--automorphism (\ie the conjugation) as
\[
\overline{l_1 + l_2Y} = \overline{l_1} - \overline{l_2}Y,\quad l_1,l_2 \in L.
\]

The conjugation appears in the definitions of the reduced trace $Tr_{B/F}$
and the reduced norm $Nr_{B/F}$ on $B$
\begin{gather*}
Tr_{B/F},~Nr_{B/F}: B \to F \\
Tr_{B/F}(b) = b + \bar b,\quad
Nr_{B/F}(b) = b\bar b,~~b \in B. 
\end{gather*}

Note that because $B$ is positive definite, the
norm of any element in it is totally positive and every $\mathfrak
o_F$--ideal in $B$ has a basis of four elements. 

The following  order will be used later on in this paragraph.  
Let
\begin{equation}\label{bigorder}
\mathcal O' = \mathfrak o_F[1,~X,\frac{w}{2}+\frac{1}{2}Y,
\frac{w}{2}X+\frac{1}{2}XY].
\end{equation}

It is easy to verify that $\mathcal O'$ is an order of $B$. We recall that
 the reduced discriminant $d_r(\Lambda)$ of an
order  $\Lambda$ in $B$, whose associated lattice has a
$\Lambda_F$-basis $\{e_1,\ldots,e_4\}$, is the integral ideal of
$\mathfrak o_F$ generated by $\sqrt{\det(Tr_{B/F}(e_i\cdot\bar e_j))}$.  
Easy computations show that the reduced discriminant of the order
$\mathcal O'$ is
\[  
d_r(\mathcal O') = \wp_2\cdot \wp_3\cdot\wp_5
\]
\ie $\mathcal O'$ has level $\wp_5$. Because $d_r$ is square free,
it follows that $\mathcal O'$ is hereditary (\cf\cite{B}, Proposition
(1.2), p.~304). 

Now, we construct a maximal $\mathfrak o_F$--sublattice
of the $\mathfrak o_F$--lattice  associated to $\mathcal O'$. Using
this lattice we plan to define a maximal Eichler suborder of $\mathcal O'$
of level $5$ (\ie with reduced discriminant $d_r = (\wp_2\cdot \wp_3)\cdot 5$). 
Consider the following  $\mathfrak
o_F$-invertible matrix  
\[
M = \begin{array}({cccc}) 1&0&0&0\\
0&1&0&0\\
0&0&2&1\\
0&0&1&w\end{array}.
\]
By letting $M$ act, by right multiplication, on the following
$4\times 4$ matrix that represents (column way) the $\mathfrak
o_F$--lattice associated to $\mathcal O'$  
\[
\begin{array}({cccc}) 1&0&w/2&0\\
0&1&0&w/2\\
0&0&1/2&1\\
0&0&1&1/2\end{array},
\]

we obtain the following suborder $\mathcal O \subset
\mathcal O'$ 
\begin{equation}\label{order}
\mathcal O = \mathfrak o_F[1,~X,~-\frac{1}{2}+\frac{3-w}{2}Y,
-\frac{w}{2}-\frac{w+1}{2}X+\frac{1}{2}Y+\frac{w}{2}XY]. 
\end{equation}

It follows from the construction that the index $[\mathfrak o:\mathfrak
o'] = \sqrt 5$. It is easy to verify that $d_r(\mathcal O) = (\wp_2\cdot
\wp_3)\cdot 5$. 

Because $B$ ramifies at $\wp_2$ and $\wp_3$, $\mathcal O$
is everywhere maximal except at the place $\wp_5$. There, the level
of $\mathcal O$ is two steps below the maximal. Finally, as a consequence of the fact
that $\mathcal O'_{\wp_5}$ is not a maximal order in a ramified
skew-field (\ie $B$ splits at $\wp_5$), we can conclude that $\mathcal
O'$ is an Eichler order. The order is isomorphic,
locally at $\wp_5$, to the order of the matrices 
\[
\begin{array}({cc}) \mathfrak o & \mathfrak o
\\ \wp_5\mathfrak o & \mathfrak o \end{array},
\]

where $\mathfrak o$ is the ring of integers of the local field
$F_{\wp_5}$. \vspace{.1in} 

It is well known that to an order $\Lambda$ of $B$ one can (locally)
associate an integer called the Eichler invariant. Let
$v$ be a finite place of $F$ and let $\mathfrak o_v$ be the completion at
$v$ of the local ring $\mathfrak o_{F,v}$, with maximal ideal $\mathfrak m$ and
residue field $F_v = \mathfrak o_{v}/\mathfrak m$. Then,  the Eichler invariant
$e(\mathcal O_{v})$ of $\mathcal O_v := B_v \otimes
\mathcal O$ ($B_v = B \otimes F_{v}$) is an
integer whose  value depends upon the description of the quotient
$\mathcal O_{v}/J_{v}$, where $J_v$ is the Jacobson radical of the
ring $\mathcal O_{v}$.  We recall from the theory of Eichler orders
developed in \cite{B}, that
\begin{equation*}
e(\mathcal O_{v}) = 
\begin{cases} 
-1& \text{if $\mathcal O_{v}/J_{v}$ is a quadratic field
extension of $(\mathfrak o_v,\mathfrak m)$}, \\  
0& \text{if $\mathcal O_{v}/J_{v} \simeq \mathfrak o_v/\mathfrak m$}, \\ 
1& \text{if $\mathcal O_{v}/J_{v} \simeq \mathfrak o_v/\mathfrak m \times
\mathfrak o_v/\mathfrak m$}.
\end{cases}
\end{equation*}

If $e(\mathcal O_{v}) = 1$, then $\mathcal O_{v}$ has exactly two
minimal overorders and it is an intersection of two uniquely
determined maximal orders. Such an order is often called Eichler order
and it is isomorphic  to the order consisting of matrices
\[
\begin{array}({cc}) \mathfrak o_v & \mathfrak o_v
\\ \pi^d\mathfrak o_v & \mathfrak o_v \end{array}
\]
where $\pi$ denotes a generator (uniformizer) of the maximal ideal $\mathfrak m$
and $d$ is a non negative integer. One can fix $\pi$ so that
$v(\pi) = 1$, where $v$ denotes the valuation of $F_v$ that corresponds to
$\mathfrak o_v$. 

From what we said, it follows that $e(\mathcal O'_{\wp_5}) = 1$. 
Hence, the value of $e(\mathcal O_{\wp_5})$ can be
either $1$ or $0$. For, if by absurd $e(\mathcal O_{\wp_5})= -1$,
then $\mathcal O'_{\wp_5}$  would necessarily  have the same Eichler
invariant (\cf~\cite{B} Prop. 3.1), in contradiction with the fact that $\mathcal O'$ is an Eichler order.  

At the two places $v = \wp_2,\wp_3$, $e(\mathcal O_v) = -1$,
as $\mathcal O_v$  is maximal there. Therefore, the
order $\mathcal O$  has level $5$. Note that the parenteses in $d_r(\mathcal O) =
(\wp_2\cdot \wp_3)\cdot 5$ separate the product of the ramified places for B, \ie
$D(B/F)$  from the level of the order. 

We continue this digression on invariants (locally) associated to
orders in quaternion algebras by recalling the definition of the mass
of an order. The mass plays a central role in the
determination of the class and the type number of an (Eichler) order.

To an order $\Lambda$ in a quaternion algebra $B$ one associates  a
rational value $m(\Lambda)$ called the mass. It is well known that the
mass  depends only upon the level (not on the particular order
chosen).  Its value appears  in the description of the class number
$h(\Lambda)$ of the order as the number of distinct classes of
left-$\Lambda$-ideals (\cf~\cite{V}). The computation of the mass
depends upon certain invariant values associated to the field
extension $F/\Q$ and upon the Eichler symbols of the order $\Lambda$
at each  prime dividing the level. A general description of the
mass  can be found in \cite{Ko} (\cf~Thm 1)
\begin{equation}\label{mass}
m(\Lambda) =
\frac{2D_K^{3/2}h_K\zeta_K(2)}{(2\pi)^{2n}}Nr_{B/F}(d_r(\Lambda))\prod_{p\vert
d_r}\frac{1-Nr(p)^{-2}}{1-e(\Lambda_p)Nr(p)^{-1}}.   
\end{equation}

Here $K$ denotes a totally real algebraic number field, $D_K$ is its
absolute discriminant, $h_K$ is the class number, $n$ is the degree of the
field extension and $\zeta_K$ is the Dedekind zeta function of $K$.

The formula \eqref{mass}, applied to our case yields, for $n=2$
and
\[
\zeta_F(\wp_2) = D_F^{-3/2}(-2\pi^2)^n \vert\zeta_F(-1)\vert,\quad
\vert\zeta_F(-1)\vert  = 1/30,\quad D_F = 5 
\]
to 
\[
m(\mathcal O) =
\frac{2(5)^{3/2}\zeta_F(-1)(5)^{-3/2}(-2\pi^2)^2}{(2\pi)^4}2^23^25^2\frac{1-2^4}{1+2^2}\cdot
\frac{1-3^4}{1+3^2}\cdot \frac{1-5^2}{1-e(\mathcal O_{\wp_5})5^{-1}} =
\]
\[
= \frac{48}{5-e(\mathcal O_{\wp_5})}.
\]
Notice that in the Eichler case (\ie for $e(\mathcal O_{\wp_5}) =
1$), one obtains $m(\mathcal O) = 12$.

The description of the class number $h(\Lambda)$ of an order $\Lambda$ 
depends upon the value of the mass. For the determination of it,
we refer to \op (\cf~Thm. 2). Before
to state the formula (in the form that
applies to our case), we introduce (following  \op) few further
notations.   

Let $B$ be the quaternion algebra described in
Lemma~\ref{quat} and let $\mathcal O$ be the order defined in
\eqref{order}. We denote by $C(\mathfrak o_F,1)$ the finite set
consisting of all $\alpha \in \mathfrak o_F$, with $\alpha^2 - 4 \notin
F_{v}^2$, where $v$ ranges among the set $S$ of places of $F$ where
$B$ is a skew-field (\ie it ramifies). It follows that $C(\mathfrak o_F,1) =
\{\frac{1 \pm\sqrt 5}{2},~\frac{-1 \pm\sqrt 5}{2}\}$. For each
$\alpha \in C(\mathfrak o_F,1)$, we consider the (separable)
quadratic extension $L = F[x]$ of $F$ defined by the quadratic equation
\[
x^2 - \alpha x + 1 = 0.
\]

For each extension $F[x]$, we define the maximal order $\Omega = \mathfrak
o_F + \mathfrak o_F[x]$. The four orders $\Omega$'s so obtained are all
isomorphic to
\[
\Omega = \mathfrak o_F + \mathfrak o_F[\frac{1}{2}(w + i\sqrt{2+w})].
\]
The extensions $L$ defined by each of the above quadratic equations
(as $\alpha$ varies in $C(\mathfrak o_F,1)$) are
unramified at the places $\wp_2,\wp_3$ and they are ramified at $\wp_5$. Under this
setting, the Theorem~2 in \op computes the class number $h(\mathcal
O)$, as a trace of the Brandt matrix $B(\mathcal O,\mathfrak o_F)$. In our case the
formula reads as
\begin{equation}\label{classn}
h(\mathcal O) = m(\mathcal O) +  \sum_{\alpha,\Omega}\prod_{v \in
S\atop v \vert d_r(\mathcal O)}
E(\Omega_v,\mathcal O_v)\frac{h(\Omega)}{2[\Omega^{\ast}:\mathfrak
o_F^{\ast}]}. 
\end{equation}

Here, $E(\Omega_v,\mathcal O_v)$ denotes  the (local)
embedding numbers (\cf~\cite{B}): \ie the number of classes of all
optimal embeddings $\Omega_v \hookrightarrow \mathcal O_v$ at each
place $v \in S$  (\cf~\cite{B2} and \cite{V}). 

To compute the embedding numbers at the two ramified places $\wp_2,\wp_3$,
one may  apply the theory resumed in Vigneras book (\cf~\op Th\'eor\`eme 3.1)
and  due to Hijikata. Because the quadratic extension $L$ is
unramified at $v = \wp_2,\wp_3$ and because the orders $\Omega_v$ and $\mathcal O_v$
are both maximal there, it follows that the number of optimal
embeddings modulo the action of the units $\mathcal O_v^{\ast}$ is
$E(\Omega_v,\mathcal O_v) = 2$. At the place
$\wp_5$, where $L$ ramifies, one can apply the theory
developed  by Brzezinski (\cf~\cite{B2} Corollary 1.6 for the Eichler
case or \cite{B1} Theorem 3.10)  to conclude that for $e(\mathcal
O_{\wp_5}) = 1$, $E(\Omega_{\wp_5},\mathcal O_{\wp_5})
= 0$.  
From what we said we conclude that
\[
h(\mathcal O) = m(\mathcal O) = 12 \qquad \text{for}\quad e(\mathcal
O_{\wp_5}) = 1. 
\]

It remains to show that $e(\mathcal O_{\wp_5}) = 1$. To this purpose,  we use the well known correspondence that relates lattices
on quadratic spaces and orders (\cf\cite{B1}). Namely,
ternary lattices define orders in quaternion
algebras. To a quaternion algebra $B$ over the field
of quotients $F$ of a Dedekind ring $R$,  
one can associate the half-regular quadratic space  $(B_0,Nr_{B/F})$, where
\[
B_0 := \{b \in B : Tr_{B/F}(b) = 0\}.
\]
Every $R$-order $\Lambda$ in $B$ defines a $R$-lattice
$\Lambda^* \cap B_0$ on $B_0$, where
\[
\Lambda^* := \{b \in B : Tr_{B/F}(b\cdot\Lambda) \subset R \}.
\]
$\Lambda^*$ is called the dual of the lattice $\Lambda$ with respect
to the trace form on $B_0$. 

Let suppose from now that $R$ is a principal ideal ring and let $e_i$
be a (finite) $R$-basis for the lattice $\Lambda^*$. To $\Lambda^*$
one associates the following ternary quadratic form
\begin{equation}\label{ternary}
q_{\Lambda^*}(X_1,\ldots,X_i,\ldots) =
\frac{1}{Nr_{B/F}(\Lambda^*)}(\sum_iNr_{B/F}(e_i)X_i^2+\sum_{i<j}
Tr_{B/F}(e_i,\bar e_j)X_iX_j). 
\end{equation}
Note that $q_{\Lambda^*}$ is (locally) $R$-integral and
it is well defined as $R$ is a principal ideal ring. Furthermore,
$q_{\Lambda^*}$ is primitive, that is, the $R$-ideal generated by its
coefficients is equal to $R$. It is well known (\cf\op) that if $p$ is
a prime ideal in $R$, then the equality $e(R_p) = 1$ implies that $q_{\Lambda^*}$ splits
in a product  of two different linear factors when its coefficients
are  evaluated in the residue field at $p$. This is what we are going
to check  in our case. 

Because the order $\mathcal O$ defined in \eqref{order} is a
Gorenstein, Bass
order (\ie $d_r(\mathcal O)$ is locally cube-free),
it follows that the level $l(\mathcal O_{\wp_5}) :=
Nr_{B/F}(\mathcal O_{\wp_5}^*)^{-1}$ is equal to the inverse of the
(local) reduced discriminant $d_r(\mathcal O_{\wp_5}) = 5$.  

We re-consider the basis chosen for the order $\mathcal O$, \ie
\begin{gather}\label{basis}
f_1 = 1,\quad f_2 = X,\\
f_3 = w + \frac{w}{2}X + Y + \frac{1}{2}XY,\quad f_4 =
\frac{w}{2}+\frac{1}{2}(w+1)X + \frac{1}{2}Y + \frac{w}{2}XY. \notag
\end{gather}

The dual of this basis, with respect to the reduced trace form (\ie
$g_i$ such that $Tr_{B/F}(f_i\cdot\bar g_j) = \delta_{ij}$) is
\[
g_1 = \frac{1}{2} - \frac{1}{10}(1 + 3w)Y,\quad g_2 = \frac{1}{12}X - \frac{1}{60}(1 + 3w)XY,
\]
\[
g_3 = \frac{1}{5}(1 + w)Y - \frac{1}{30}wXY,\quad g_4 =
-\frac{1}{15}w(3Y - XY).
\]

We represent this basis by means of a $4\times 4$ matrix whose 
i-th column stands for $g_i$ written down on the basis $X$, $Y$, $XY$
and $1$. After few column reductions, this matrix takes the following
upper triangular form 
\[
\begin{array}({cccc}) \frac{1}{6}(2-w)&0&\frac{1}{12}(-3+2w)&0\\
0&-\frac{1}{5}(1+3w)&\frac{1}{5}(1+w)&\frac{1}{10}(1+3w)\\
0&0&\frac{1}{60}(1+w)&0\\
0&0&0&\frac{1}{2}\end{array}.
\]
The first three columns represent a basis for the $\mathfrak
o_F$-lattice $\mathcal O^* \cap B_0$. As elements of the quaternion
algebra, they are described, in terms of the chosen quaternion basis, as
\[
h_1 = \frac{1}{6}(2-w)X,\quad h_2 = -\frac{1}{5}(1+3w)Y,
\]
\[
h_3 = \frac{1}{12}(-3+2w)X + \frac{1}{5}(1+w)Y + \frac{1}{60}(1+w)XY.
\]
The corresponding ternary quadratic form is
\begin{multline*}
q_{\mathcal O^*\cap B_0} = \frac{5}{6}(5-3w)X_1^2 + (3+4w)X_2^2 + \frac{5}{6}(4-w)X_3^2 + \frac{5}{6}(-8+5w)X_1X_3 -2(1+2w)X_2X_3.
\end{multline*}
The evaluation of this form in the residue field at $\wp_5$ gives
\[
\bar q_{\mathcal O^*\cap B_0} = X_2X_3.
\]
This shows that $e(\mathcal O_{\wp_5}) = 1$. Hence, we can conclude that  $\mathcal O$ is an Eichler order.

Finally, using Theorem~ 3 in \cite{Ko}, we  
determine the type number $t(\mathcal O)$ of $\mathcal O$. 
This quantity counts the number of non equivalent types of Eichler
orders with reduced discriminant $(\wp_2\cdot \wp_3)\cdot 5$ in
the chosen quaternion algebra $B$. 
The formula for the computation of the type number reads in our case
as ($S = \{\wp_2,\wp_3,\wp_5\}$) 
\begin{equation}\label{type}
m(\mathcal O) = t(\mathcal O)\prod_{v\in S\atop v \vert d_r(\mathcal
O)}[\Gamma(\mathcal O_v):F_v^*\mathcal O_v^*]. 
\end{equation}
Here, $[\Gamma(\mathcal O_v):F_v^*\mathcal O_v^*]$ is the cardinality
of the central Picard group of $\mathcal O_v$, where $v$ denotes a
place in $S$.  It is known (\cf~\cite{Ko1} Satz 2)
that this number is $1$ when $\mathcal O$ is maximal at the place
$v$. On the other hand, for orders whose Eichler invariant is non zero
at a place $v$, the cardinality of the central Picard group is either
$1$ or $2$, depending upon the value of $ord_v(d_r(\mathcal O_v))$
(\ie being $0$ or not). Easy computations show that in our case
$t(\mathcal O) = 3$. 

We end this paragraph by resuming the main results proved it it
\begin{prop} Let $B$ be the quaternion algebra with discriminant $D(B/F)
= \wp_2\cdot \wp_3$ defined in
Lemma~\ref{quat} and let $\mathcal O$ be the Eichler order of level
$5$, defined in \eqref{order}. Then
\[
h(\mathcal O) = 12,\quad t(\mathcal O) = 3.
\]
\end{prop}

\section{Hilbert modular forms and quaternion algebras.}\label{4} 


In this paragraph we review the description
of the space
of quaternionic cusp forms, namely functions defined on a quaternion
algebra. One of the main results in the theory is the
Jacquet-Langlands Theorem. Shorthly said, this theorem relates the systems of
eigenvalues found in the spaces of quaternionic cusp forms to certain
systems of eigenvalues 
attached to the spaces of classical modular forms. In this paper, the term
classical modular forms
means holomorphic Hilbert cusp forms over $F =
\Q(\sqrt 5)$ of weight ${\bf k} = (2,4)$ and level 
$\mathfrak N := \wp_2\wp_3\cdot 5$: an integral ideal of $\mathfrak o_F$.  
Our exposition on this correspondence follows the report 
given by R.~Taylor in \cite{T} and by Hijikata, Pizer and
Shemanske in \cite{HPS}. For self-contained reasons
we restrict the description to the case of
interest. 

First, we introduce some notations to which we refer
throughout this paragraph and later on in the paper. 

Let $F$ denote the totally real quadratic field $\Q(\sqrt 5)$ and let
$I = \{\iota_\tau\}_{\tau = 1,2}$ be
the set of the embeddings $F \hookrightarrow \R$. We shall let $B$
denote the definite quaternion algebra with center $F$ 
and discriminant $D(B/F) = \wp_2\wp_3$ introduced in Lemma~\ref{quat}. 
We write $\mathcal O$ for the Eichler order of $B$ of level $5$ defined in
\eqref{order}.   

We will consider the following algebraic groups: $\mathfrak G = {\bf GL}_2$
and $G_B$ over 
$F$: $\mathfrak G(F) = {\rm GL}_2(F)$ and $G_B(F) = B^\times$. 
Their reduced norm morphisms are $\nu: \mathfrak G \to {\bf G}_m$ and $\nu_B:
G_B \to {\bf G}_m$. Here ${\bf G}_m$ denotes the multiplicative group
in one variable considered as an algebraic group defined over $F$.
 
Let $\A_F =  \A_f \times \A_\infty$ denote the ring of adeles of $F$,
decomposed into its  finite and infinite parts. For an algebraic group
${\bf G}$ over $F$, we write ${\bf G}_f = {\bf G}(\A_f)$ and ${\bf
G}_{\infty} = {\bf G}(\A_\infty)$. 

For $v$ a finite place of $ F$ such that $v \nshortmid
d_r(\mathcal O)$, we fix isomorphisms $M_2(\mathfrak o_{F,v}) \simeq
\mathcal O \otimes_{\mathfrak o_F} \mathfrak o_{F,v}$. They induce
the isomorphism $(G_B)_f \simeq \mathfrak G_f$. We choose a subfield $K$ of $\C$ 
which is Galois over $\Q$ and that splits $B$, so that there is an
isomorphism $i: \mathcal O \otimes_{\Z} \mathfrak o_K \to M_2(\mathfrak o_K)^I$. 

In what it follows we fix ${\bf k} = (2,4) \in \Z^I$.
\vspace{.1in}

First, we describe a model  for the space of Hilbert modular forms of
weight ${\bf k}$. This is done by illustrating the action of the Hecke
algebra on it. Our description is given in terms of automorphic
functions on the adele group of $\mathfrak G$ and the corresponding
Hecke algebra.

Let us consider functions as $f: \mathfrak G(\A_F) \to \C$. For $u =
u_fu_\infty \in \mathfrak G_f \times \mathfrak G_\infty$, define
\[
(f_{\vert_{\bf k}})(x) =
j(u_\infty,z_0)^{-k}\nu(u_\infty)^{\frac{k}{2}}f(xu^{-1}),
\]
where
\begin{align*}
&a)\quad z_0 = (\sqrt{-1},\sqrt{-1}) \in \mathcal Z^I,~\mathcal
Z~\text{upper half complex plane} \\
&b)\quad j: \mathfrak G_\infty \times \mathcal Z^I \to \C^I,~\left (\begin{array}{cc}
a_\tau & b_\tau \\
c_\tau & d_\tau 
\end{array}\right ) \times z_\tau \mapsto (c_\tau z_\tau + d_\tau).
\end{align*}

For $U \subset \mathfrak (G_B)_f$ an open compact subgroup, one defines
\[
\mathcal S_{{\bf k}}(U) := \{f: \mathfrak G(F)\backslash \mathfrak G(\A_F) \to \C~\vert~\text{under the
following conditions 1.,2.,3.}\}
\]
\begin{equation*}
\begin{split}
1.&~f_{\vert_{{\bf k}}}u = f~\text{for all}~u \in
  UC_\infty,~\text{where}~C_\infty = (\R^\times\cdot SO_2(\R))^I
  \subset \mathfrak G_\infty,    \\
2.&~\forall x \in \mathfrak G_f,~f_x:\mathcal Z^I \to
  \C:~uz_0 \mapsto j(u,z_0)^k\nu(u)^{-\frac{k}{2}}f(xu),~u \in
  \mathfrak G_\infty,~\text{is holomorphic,} \\
3.&~\int_{\A_F\slash F}f\left(\left (\begin{array}{cc}
1 & a \\
0 & 1 
\end{array}\right )x\right )da = 0,~\forall~x \in
\mathfrak G(\A_F)~da~\text{additive Haar measure on}~\A_F\slash F.
\end{split}
\end{equation*}

It is well known that the space $\mathcal S_{\bf k}(U)$ is a model for 
the the space of holomorphic Hilbert modular
cusp forms on $F$ of weight ${\bf k}$. The level of these forms is directly related to 
a choice of the open set $U$. We refer to \cite{VdG} (Definition~$6.1$) for a definition of the
space of Hilbert modular forms. 
In this model, the action of $\mathfrak G_f$ is given by
right translation.

We define $U_0 = \prod_\wp \mathfrak G(\mathfrak o_{F,\wp})$, where
$\wp$ runs over the set of finite, integral, prime ideals in $F$. From
now on we write 
$\mathfrak P$ for the ideal $\wp_2\wp_3$ in $\mathfrak o_F$ product of the two
prime inert ideals above $2$ and $3$ in $\Z$; we set $\mathfrak N = \mathfrak
P\cdot 5$. Consider the
following open compact subgroup of $(G_B)_f$
\[
U(\mathfrak P,5) = \{\left (\begin{array}{cc}
a & b \\
c & d 
\end{array}\right ) \in U_0~\vert~c \in \mathfrak P\cdot 5,~a-1 \in
\mathfrak N\}.
\]

We will denote by $\mathcal S_{\bf k}(\mathfrak P,5) = \mathcal S_{\bf k}(U(\mathfrak
P,5))$. The space $\mathcal S_{\bf k}(\mathfrak P,5)$ 
is in one-to-one correspondence with the $\C$--vector space of holomorphic
Hilbert modular cusp forms on $F$ of weight ${\bf k}$ and level
$\mathfrak N$.  

The description of the Hecke algebra on the space $\mathcal S_{\bf k}(\mathfrak P,5)$ 
goes as follows.

For a prime ideal $\wp$ of $F$ and
for the choice of the open compact $U = U(\mathfrak P,5)$, set
\[
T_\wp = \left [U\left (\begin{array}{cc}
1 & 0 \\
0 & \pi_\wp 
\end{array}\right )U \right ].
\]
$T_\wp$ denotes the Hecke operator at $\wp$, where $\pi_\wp \in \A_f$ is a
uniformizer at $\wp$. There is a similar definition
of the Hecke operator, when $\wp$ is replaced by a fractional ideal
$\mathfrak a$ of $F$ with $(\mathfrak a,\mathfrak N) = 1$. Note that, although
the uniformizer $\pi_\wp$ is not uniquely 
determined by $\wp$, $T_\wp$ is well defined for the choice we made of
the open set $U$ (a similar remark holds for the Hecke operator
associated to a fractional ideal $\mathfrak a$ of $F$).

As $\wp$ ranges among all finite prime ideals
of $F$ (and 
among all finite integral ideals of $F$ prime to $\mathfrak N$
), the $\Z$--algebra ${\bf T}_{\bf k}(\mathfrak N) \subset {\rm End}(\mathcal S_{\bf k}(\mathfrak
P,5))$ of all 
these Hecke operators  is known to be diagonalizable on the (unique) ${\bf T}_{\bf k}(\mathfrak
N)$ submodule of $\mathcal S_{\bf k}(\mathfrak P,5)$ generated by the eigenforms of
the Hecke algebra (\ie $f \in \mathcal S_{\bf k}(\mathfrak P,5)$, such
that $f_{\vert{\bf k}}T = \theta_f(T)f$, $\forall$ $T \in {\bf T}_{\bf k}(\mathfrak N)$). 

The space $\mathcal S_{{\bf k}}(\mathfrak P,5)$ may be equivalently described via the
introduction  of automorphic forms
defined on $G_B(\A_F)$. We will recall this construction, 
together with the corresponding algebra of Hecke operators and few of
their properties. 
  
We denote by $S_{a,b}(\C)$ the right
${\rm M}_2(\C)$-module ${\rm Sym}^a(\C^2)$ endowed with the following
${\rm M}_2(\C)$--action ($m \in {\rm M}_2(\C)$, $s \in {\rm Sym}^a(\C^2)$)
 
\[
s\cdot m : = (\det m)^b~s~{\rm Sym}^a(m)
\]
For ${\bf k} = (2,4)$, we set
\begin{equation}\label{kk}
L_{{\bf k}}(\C) = S_{0,1}(\C) \otimes S_{2,0}(\C).
\end{equation}

Each factor of the tensor product is associated to a
real embeddings $\iota_\tau \in I$. The group ${\rm M}_2(\C)$ acts (on the
right) on $S_{0,1}(\C)$ 
by means of the determinant and via the second symmetric square on
$S_{2,0}(\C)$. Hence, the space $L_{\bf k}(\C)$ inherits a
$(G_B)_\infty$ action via 
the injection $i: (G_B)_\infty \hookrightarrow
\mathfrak G(\C)^I $ that we have fixed at the beginning of this paragraph. 

Now, we consider functions $f: G_B(\A_F) \to L_{{\bf k}}(\C)$. For $\tilde u =
u_fu_\infty \in G_B(\A_F)$, 
we define 
\[
(f_{\vert_{{\bf k}}}\tilde u)(x) = f(x\tilde u^{-1})\cdot u_\infty.
\]

The following open compact subgroup of $(G_B)_f$ 
``corresponds'' to $U(\mathfrak P,5)$
\[
\mathcal U = \prod_{\wp \in F\atop \wp < \infty}(\mathcal
O \otimes_{\mathfrak o_F}\mathfrak o_{F,\wp})^{\times} = \prod_{\wp \in
F\atop \wp < \infty}\mathcal O_{\wp}^{\times} = \{\tilde u =
(u_v) \in \prod_v B_v^\times~\vert~u_v \in \mathcal O_v^\times, \forall v < \infty\}.
\]

Note that the local description of the order $\mathcal O$ is
(\cf~section~\ref{4})  
\[
\mathcal O_{\wp}^\times \simeq \begin{cases} 
\mathfrak G(\mathfrak o_{F,\wp})& \text{for $\wp \nmid 5$}, \\\\  
\left \{\left [\begin{array}{cc}
a & b \\
5c & d 
\end{array}\right ] \in \mathfrak G(\mathfrak o_{F,\wp_5})~\vert~a,b,c,d \in \mathfrak
o_{F,\wp_5}\right \}& \text{at $\wp = \wp_5$}.
\end{cases}
\] 

We set
\begin{multline}\label{hmf}
S_{{\bf k}}^B(\mathcal U) := \{f: B^\times\backslash G_B(\A_F) \to
L_{{\bf k}}(\C)~\vert~f_{\vert_{{\bf k}}}\tilde u = f,~\forall
\tilde u \in \mathcal U(G_B)_\infty\} =  
\\
= \{f: (G_B)_f/ \mathcal U \to L_{\bf k}(\C),~\vert~f(b\tilde x) = f(\tilde
x)\cdot b^{-1},~\forall \tilde x \in (G_B)_f,~b \in B^\times\}.
\end{multline}

The global units $B^\times$ act via the map $i: B^\times \hookrightarrow
\mathfrak G(K)^I$ that we have chosen at the beginning of this section. 
It is worthwhile to remark that the group of the ideles of $B$ can be
written as a finite union of distinct double cosets of $\mathcal U$
and $D^\times$ ($h(\mathcal O) = 12$) 
\[
(G_B)^\times =  \coprod_{\lambda =
1}^{12}\mathcal U \tilde g_\lambda B^\times,
\]

where the representatives $\tilde g_\lambda = (g_{\lambda_\wp})$ can
be (and are) chosen so that $g_{\lambda_\wp} \in \mathcal
O^\times_\wp$ for all $\wp \vert \mathfrak N$. As $\lambda$ runs between $1$ and $12$,
let $I_\lambda = \mathcal O\tilde g_\lambda$. The double coset
space $\mathcal  X(\mathcal U) =  B^{\times}\backslash  (G_B)_f/
\mathcal U$ can be canonically identified 
 with the right equivalence classes  of left $\mathcal O$--ideals
 $I_\lambda$ and in
 particular it is finite. The elements $\tilde g_\lambda \in \mathcal
 X(\mathcal  U)$ correspond to the left ideals $I_\lambda$ of
$\mathcal O$, each of  which determines a different class in the
order. The elements of $S^B_{{\bf k}}(\mathcal U)$ are completly
determined by their values at $\tilde g_\lambda$. In other words, the
setting $f \mapsto (f(g_\lambda))$  defines the isomorphism  
\begin{equation}\label{ident}
S_{{\bf k}}^B(\mathcal U) \rightarrow 
\bigoplus_{\lambda=1}^{12}(L_{{\bf k}}(\C))^{\mathcal U~\cap~\tilde
g_\lambda^{-1}B^\times\tilde g_\lambda}.
\end{equation}

In analogy with the previous construction, one may define the Hecke algebra
${\bf T}_{\bf k}^B(\mathcal U) \subset {\rm End}(S_{\bf k}^B(\mathcal U))$ as the
algebra generated by the Hecke operators
\[
\left [\mathcal U\tilde g\mathcal U'\right ]: S_{\bf k}^B(\mathcal U) \to
S_{\bf k}^B(\mathcal U'),\quad f \mapsto \sum
f_{\vert_{\bf k}}\tilde g,
\]
for $\mathcal U'$ open compact subgroup and
$\mathcal U\tilde g\mathcal U' = \coprod_i \mathcal U gg_i$ the decomposition
in double cosets. In particular, for $\mathcal U = \mathcal U'$, we
have  $[\mathcal U\tilde g\mathcal U] \in {\bf
T}_{\bf k}^B(\mathcal U)$, where $\mathcal U\tilde g\mathcal U = \coprod_i\mathcal
Ug_i$ is the decomposition into disjoint right
cosets. The representation $\rho$ of ${\bf T}_{\bf k}^B(\mathcal U)$ on
$S_{\bf k}^B(\mathcal U)$ is  defined as follows
\[
\rho([\mathcal U\tilde g\mathcal U])f(x) =  (f_{\vert_{\bf k}}\tilde g)(x) =
\sum_i f(g_i x), \qquad \tilde g = (g_i).
\]
Then, one  extends $\rho$ to all of ${\bf T}_{\bf k}^B(\mathcal U)$ by linearity.
Because of the identification \eqref{ident}, we set $f_\lambda = f(\tilde g_\lambda)$,
for  $f \in S_{\bf k}^B(\mathcal U)$.
Then, the mapping $\alpha: f \mapsto
(f_1,\ldots,f_{12})$ describes an isomorphism of $S_{\bf k}^B(\mathcal U)$ into
$\C^3 \oplus \C^3 \oplus\ldots\oplus \C^3$. This sum has exactly $h(\mathcal O)
= 12$ addenda. By using this identification, one can think of $\rho$
as giving a matrix representation of ${\bf T}_{\bf k}^B(\mathcal U)$ on
$(\C^3)^{12}$. 
Let $\xi \in \mathfrak o_{F,>>}$ be a
totally positive element in the ring of integers of $F$. 
The corresponding Hecke operator is defined as
\[
T(\xi) = \sum_{Nr(\tilde g) \in \xi \mathcal U}[\mathcal U\tilde g
\mathcal U],\quad \text{for $\tilde g
= (g_\wp),~g_{\wp} \in \mathcal O_\wp,~\forall \wp < \infty$}.
\]

Here $Nr$ denotes the usual extension of the norm $Nr_{B/F}$ to $G_B(\A_F)$.
One may represent $\rho(\xi) = \rho(T(\xi))$ as the $12\times
12$ (block) matrix 
\begin{align}\label{br}
&B(\xi) = [\rho_{i,j}(\xi)]_{i,j=1,\ldots 12},~\rho_{i,j} =
\text{pr}_i\cdot\alpha\cdot\rho(\xi)\cdot\alpha^{-1} \\ 
&\rho_{i,j}: \C^3_j \hookrightarrow (\C^3)^{12} \to S_{\bf k}^B(\mathfrak P,5) \to
End(S_{\bf k}^B(\mathfrak P,5)) \to (\C^3)^{12} \twoheadrightarrow \C^3_i.\notag
\end{align} 
The first map on the left denotes the injection of the j-th factor of
$(\C^3)^{12}$ into the whole space. The following maps are \resp
$\alpha^{-1}$, $\rho(\xi)$, $\alpha$ and finally the i-th projection of
$(\C^3)^{12}$ onto $\C^3$.
In the next paragraph we will give an explicit description of  the
matrix $B(\xi)$.

We refer to \cite{G} for a proof of the (Hecke equivariant)
isomorphism $S_{\bf k}^B(\mathcal U) \simeq \mathcal S_{\bf k}(\mathfrak
P,5)$. Here, we simply recall the important fact (for which we 
refer again to \op) that these spaces
have  an 
equivalent description in terms of representations of
\resp the groups $G_B(\A_F)$ and $\mathfrak G(\A_F)$. To the ``new-forms'' in the
spaces of automorphic forms (\ie eigenfunctions for the Hecke
operators) correspond irreducible
representations of the groups of the adeles of $G_B$ and $\mathfrak G$.

The second part of this paragraph is expository in nature. We plan to describe how the Jacquet-Langlands
construction may help in the process of searching for a particular Hilbert
modular form. Roughly said, the idea is to ``translate'' a local information supplied
by the Galois representation $\rho$ (\cf paragraph~\ref{J3}) at each place $v$ where $\rho$ ramifies, to  an information regarding the properties of the local component $\pi_v$ of an
automorphic representation $\pi = \otimes_v \pi_v$ of $G_B(\A_F)$ 
and  the properties of the underlying quaternion
algebra and order.  Note that this process makes sense and it 
is motivated by the knowledge that 
$\rho \otimes \Q_\ell(\sqrt 5) = \Ind_{F}^{\Q}(\sigma)$ (\cf theorem~\ref{Jaspth}). 

On the Galois-side, important
informations are encoded in the description of the Gamma-factors and  
in the (Artin) conductor exponents of the
functional equation attached to  
$\rho$. These informations  determine (uniquely) the weight ${\bf k}$ of a modular form over $F$,
the ramification places of a quaternion algebra $B$ and
the level of an associated Eichler order $\mathcal O$.

The Jacquet--Langlands Theorem  establishes
a (Hecke equivariant) injection $\pi \mapsto JL(\pi)$ from    
the set of (classes of) automorphic representations $\pi = \otimes_v
\pi_v$ of $G_B(\A_F) = (B
\otimes_F \A_F)^\times$ with $\dim \pi > 1$ (\ie not
characters), to the set of
cuspidal automorphic representations of
$\mathfrak G(\A_F)$. The theorem characterizes the image
of the ``JL'' map as the set 
of cuspidal automorphic representations of $\mathfrak G(\A_F)$ which
are discrete series (\ie special or supercuspidal at a finite
place) at all places at which 
$B$ ramifies. 
The representation $JL(\pi)$ is locally
characterized, at a place $v$, only in terms of the related $\pi_v$, in the
sense that the image only depends on $\pi_v$. If $v$ is a finite place
at which $B$ splits, then  
$JL(\pi)_v = \pi_v$. 

The description of $\pi_\infty$ is a consequence of
the following remarks. 
As explained earlier on in the paper, the expected modularity of the two-dimensional representation
$\sigma$ of
$Gal(\bar F/F)$ is equivalent to the existence of a
holomorphic Hilbert modular form $\mathfrak f$ of weight ${\bf k} =
(2,4)$. This particular choice for the weight is suggested  by the
shape of the gamma--factors in the functional equation of the
$L$--function  associated to
$\rho$ (\cf section~\ref{J4}) and  by the fact that 
the functional equations of $\rho$ and $\sigma$ coincide.  
A holomorphic Hilbert modular form of 
weight ${\bf k}$ on $F$ verifies the following functional equation  (\cf~\cite{Sh})  
\[
L(\mathfrak f,s) = L(\mathfrak f,4-s)
\]
where
\[
L(\mathfrak f,s) : = Nr_{F/\Q}(\mathfrak
c\delta^2)^{s/2}(2\pi)^{-2s}\Gamma(s)\Gamma(s-1) D(\mathfrak f,s).
\]

We denoted by $\mathfrak c$ the conductor of the (Galois) representation
associated to $\mathfrak f$ (\ie the level
$N$ of $\mathfrak f$) and by $\delta = \wp_5 = \sqrt
5$ the  different of the extension $F/\Q$. The symbol $D(\mathfrak f,s)$ means the
product of the  local, non Archimedean Euler factors. 
As we explained in details in paragraph~\ref{J4}, our numerical tests determined 
the conductor of $\rho$ to be $2^2\cdot 3^2\cdot 5^4$. Because the scalars extension
$\rho \otimes \Q_\ell(\sqrt 5)$ is induced, 
the level of $\mathfrak f$ is forced to be $\mathfrak c  = (\wp_2\cdot
\wp_3)\cdot 5$. As a consequence of the chosen weight ${\bf k}$,
the Archimedean component $\pi_{\infty}$ turns out to 
be isomorphic to the representation
\[
\sigma_{{\bf k}}': (G_B)_\infty \hookrightarrow \mathfrak G(\C)^I \to
Aut(L_{{\bf k}}(\C)),
\]

where the $\C$--module $L_{{\bf k}}(\C)$ was defined in \eqref{kk}.
By
applying the Jacquet--Langlands map, one concludes that
$JL(\pi)_{\infty} \simeq
\sigma_{{\bf k}}$ is the discrete series representation of
$\mathfrak G(\R)^I$ whose Langlands' parameter at each archimedean place
is (${\bf k} = (k_\tau) = (2,4)  \in \Z^I$)
\begin{align}\label{Lpar}
W_{\R} &= <\C^\times, j~\vert~j^2 = -1; jzj^{-1} = \bar z~\forall z \in
\C^\times >~\to \mathfrak G(\C)^I \\\notag \\
z &\mapsto  \left (\begin{array}{cc}
z^{1-k_\tau} & 0 \\
0 & \bar z^{1-k_\tau} 
\end{array}\right )(z\bar z)^{c_\tau}\vert z\vert,\quad 
j \mapsto \left (\begin{array}{cc}
0 & 1 \\
(-1)^{1-k_\tau} & 0 
\end{array}\right ).\notag 
\end{align}

Here, $(c_\tau) = (1,0) \in \Z^I$ and $\bar{}$ denotes the complex
conjugation.

At the primes $p = 2,3$ our numerical computations determined the
shape of the local Euler factors of the 
$L$-function $L(H^3(\tilde X),s)$ to be of type
$(1-p^{2-2s})^{-1}$ (\ie trivial central character). This implies
that $\pi_v \simeq \det$. Hence, $JL(\pi)_v$ is the special
representation which is a subquotient of the principal series
representation associated to the pair of characters
$\vert\cdot\vert^{-1/2}$, $\vert\cdot\vert^{1/2}$. At 
the ramified place (for the extension $F/\Q$) $\wp_5$ (the prime over
$5$) the local Euler factor was evaluated (numerically) to be 1. In order to test whether
$JL(\pi_{\wp_5})$ is a principal or a discrete series (\ie special
or supercuspidal), one would need to obtain more
information on the $L$-function. We decided to twist this function with 
the characters associated to the cyclotomic
extension of $F$ with Galois group $(\Z/5\Z)^\times$ with the purpose of testing of any
change occurring in the local Euler factor.
Unfortunately, our numerical computations did not end up with good candidates for the twisted local factors. 
The problem, however, seems  to be not so serious and we are confident that it could be definitively solved with a major accuracy in the setting of the computer program that tests the numerical convergence of the $L$-series. 
In fact, twisting the $L$-function has as consequence that of a possible increase of its conductor. In these cases, more terms in the $L$-series are needed in order to test their convergency.

\section{The Brandt matrices.}\label{5} 


In  this paragraph we describe the action
of the Hecke operators on the space of forms on a
definite quaternion algebra. More precisely, we will define
certain square matrices, naturally associated to an Eichler order on a
quaternion algebra. We  refer to the ``basis problem'' theory
 for a proof of the well-known fact that these
matrices determine on $S_{\bf k}^B(\mathcal U)$ (\cf section~\ref{4} for notations) the same action as the
matrices $B(\xi)$ in \eqref{br}. The theta--series
attached to a positive-definite quadratic form define automorphic forms
associated to certain open compact subgroups of $G_f$ (\eg
$U(\mathfrak P,5)$). Roughly speaking, the ``basis problem'' asks whether
all elements of of a 
certain space of modular forms
may be expressed as linear combinations of 
theta-series attached to the norm form of a definite quaternion
algebra (\cf~\cite{HPS}). For elliptic
modular forms this problem has been studied in great details by Hecke,
Eichler, Waldspurger, Hijikata, Pizer and Shemanske. The theory for other
modular forms has been started only recently. The
works of Shimizu (\cite{Shi}) and Waldspurger (\cite{W}) give a very useful
insight for the solution of the
``basis problem'' for Hilbert 
modular forms and for modular forms over ${\bf GL}_2$. In order to  
give an explicit construction of the space of automorphic forms on which the ``JL''
correspondence acts, Shimizu
used the spherical functions (\ie functions associated to automorphic forms defined on the
group of adeles of $G_B$) and utilized the theta series defined by
A. Weil. He showed that this space consists of cusp forms.
Waldspurger used the theory of ``newforms'' (\ie ``classes'' of common
eigenfunctions for the Hecke operators in the sense of \cite{AL}) in the proof of the ``basis problem''. 
His proof can be generalized in order to give a solution of  the ``basis
problem'' for Hilbert modular forms (\cf~\cite{W}). For other type of modular 
forms (\eg Siegel modular forms) the results are less complete
although in progress (\cf~\eg~\cite{W},\cite{B}).

The Brandt matrices describe a
representation of the Hecke operators acting on a space of theta
series. They satisfy certain commutation relations analogue to those
verified by the Hecke operators. Their eigenvalues are
subjected to certain  constraints expressed by the analogue 
of the Ramanujan-Petersson conjecture.
The aim of this paragraph is to define these matrices in our contest.
The explicit realization of them will be given in the next paragraph.

We keep  the same notations introduced in the last section.
Let start by recalling  the following construction. We refer to
\cite{P} for more details.

Consider the following reductive group over $\Q$: $G =
{\rm Res}_{F/\Q}(B^{\times})$. The group $G(\R)$ is therefore isomorphic to
$(\mathbb H^{\times})^2$, where $\mathbb H$ denote
the field of Hamilton's quaternions. The $\C$--algebra $\mathbb H$
can be represented as the subalgebra of ${\rm Mat}(2,\C)$ 
\[
\mathbb H = \{\begin{array}({cc}) z & w
\\ -\bar w & \bar z 
\end{array}
\vert~z, w \in \C \}.
\]

In this representation, a basis of $\mathbb H$ is
\[
\bf 1 = \begin{array}({cc}) 1 & 0
\\ 0 & 1 
\end{array},\quad \bf I = \begin{array}({cc}) i & 0
\\ 0 & -i 
\end{array}, \quad \bf J = \begin{array}({cc}) 0 & 1
\\ -1 & 0 
\end{array},\quad \bf K = \begin{array}({cc}) 0 & i
\\ i & 0 
\end{array}.
\]
This determines a  2-dimensional (complex) representation
\[
\Phi_1: \mathbb H^\times \to {\rm GL}(V).
\]

Let choose an isomorphism $V \simeq \C^2$. In terms of the canonical
basis $e_1 = (1,0)$, $e_2 = (0,1)$ of $\C^2$, the  matrix
representation $X_1$ of $\Phi_1$is ($\alpha = x_1 + x_2\cdot$ $\bf{I}$
$+ x_3\cdot$$\bf{J}$ $+ x_4\cdot$$\bf{K}$$ \in \mathbb H^\times$, $z
= x_1+ix_2$, $w= x_3+ix_4$)
\[
X_1(\alpha) = \begin{array}({cc}) x_1 & 0
\\ 0 & x_1 \end{array} + x_2 \begin{array}({cc}) i & 0
\\ 0 & -i \end{array} + x_3 \begin{array}({cc}) 0 & 1
\\ -1 & 0 \end{array} + x_4 \begin{array}({cc}) 0 & i
\\ i & 0 \end{array} = 
\]
\[
= \begin{array}({cc}) x_1+ix_2 & x_3+ix_4
\\ -x_3+ix_4 & x_1-ix_2 \end{array} = \begin{array}({cc}) z & w
\\ -\bar w & \bar z \end{array}.
\]

The representation $\Phi_1$ induces the determinant representation
(character)
\[
det: \mathbb H^\times \to \C^\times.
\]
Its description, always in terms of the choice of the
canonical basis $\{e_1,e_2\}$ of $\C^2$, is
\[
X(\alpha) := {\rm det}(X_1)(\alpha) = x_1^2 + x_2^2 + x_3^2 + x_4^2.
\]

Furthermore, $\Phi_1$ induces a
$3$-dimensional representation of $\mathbb H^\times$ on the 2-th
symmetric power ${\rm Sym}^2(\C^2) = (\C^2 \otimes\C^2)/K$ ($K$ is the symmetric kernel)
\[
\Phi_2: \mathbb K^\times \to {\rm GL}(\C^{3}).
\]
A basis for the space ${\rm Sym}^2(\C^2)$  is given by the set of elements
$\{e_1^ie_2^{2-i} \vert i = 0,\ldots 2\}$ (the product is
considered in the symmetric algebra). The corresponding matrix
representation is (the following formula should be read mod $K$)
\[
X_2(\alpha)(e_1^i\otimes e_2^{2-i}) = \bigotimes_1^i(X_1(\alpha)e_1)
\otimes\bigotimes_1^{2-i}(X_1(\alpha)e_2). 
\]

In the previous section (\cf~\eqref{kk}) we  introduced the notation
$L_{\bf k}(\C)$ (${\bf k} = (2,4)$) to indicate the following tensor product of (complex)
representations (each factor corresponds to an embedding of
$\Q(\sqrt 5)$ in $\C$)
\[
L_{\bf k}(\C) = S_{0,1}(\C) \otimes S_{2,0}(\C).
\] 

Here, $S_{0,1}(\C) = \C$ denotes the identity representation endowed with a
right ${\rm M}_2(\C)$--module action given by the determinant ($m \in
{\rm M}_2(\C)$, $1 \in {\rm Id}(\C)$)
\[
1\cdot m := \det m.
\]
We denote by $S_{2,0}(\C)$ the right ${\rm M}_2(\C)$--module
${\rm Sym}^2(\C^2)$. ${\rm M}_2(\C)$ acts on it on the right as
($s \in {\rm Sym}^2(\C^2)$)
 \[
s\cdot m : = s~{\rm Sym}^2(m).
\]

In the following, we  give an explicitely description of the matrix $X_2$ attached to the
representation $\Phi_2$.  

With the choice of the canonical basis $\{e_1,e_2\}$ for $\C^2$,  and
in terms of the matrix representation 
$\alpha =$ \begin{math}\left(\begin{smallmatrix}z&w\\-\bar w&\bar
z\end{smallmatrix}\right)\end{math} ($z = x_1 + ix_2$, $w = x_3 +
ix_4$, for $\alpha =  x_1 + x_2\cdot$ $\bf{I}$ $+ x_3\cdot$$\bf{J}$ $+
x_4\cdot$$\bf{K}$$ \in \mathbb H$), $\Phi_2(\alpha)(e_2^2)$ is described by 
\[
X_2(\alpha)(e_2\otimes e_2) = (X_1(\alpha)e_2) \otimes 
(X_1(\alpha)e_2) = (-\bar we_1 + \bar ze_2)^2 =
(-x_3+ix_4,x_1-ix_2)^2.
\]
Similarly, $\Phi_2(\alpha)(e_1e_2)$ is given by
\begin{multline*}
X_2(\alpha)(e_1\otimes e_2) = (X_1(\alpha)e_1) \otimes
(X_1(\alpha)e_2) = (ze_1 + we_2)\cdot(-\bar we_1 + \bar ze_2) 
= \\
= (x_1+ix_2,x_3+ix_4)\cdot(-x_3+ix_4,x_1-ix_2),
\end{multline*}

and $\Phi_2(\alpha)(e_1^2)$ takes the form
\[
X_2(\alpha)(e_1\otimes e_1) = (X_1(\alpha)e_1) \otimes 
(X_1(\alpha)e_1) = (ze_1 + we_2)^2 = (x_1+ix_2,x_3+ix_4)^2. 
\]
By patching together (row way) these three vectors, we get the
description  of the $(3\times 3)$ matrix $X_2(\alpha) = X_2(x_1 +
x_2\cdot\bf{I}$$ + x_3\cdot\bf{J}$$ + x_4\cdot\bf{K}$) as 
\begin{equation}\label{sym}
(p_{ij}(X))_{i,j=1,\ldots,3},\qquad\text{for $X = (x_1,\ldots,x_4)$.}
\end{equation}
The entries $p_{ij}(X)$ are the following quadratic, harmonic
polynomials with complex coefficients 
\begin{align*}
p_{11}(X) &= x_1^2-x_2^2+2ix_1x_2, & \quad  {p_{12}(X)} &=
2(x_1x_3-x_2x_4+i(x_1x_4+x_2x_3)),\\
p_{21}(X) &= -x_1x_3-x_2x_4+i(x_1x_4-x_2x_3),& \quad  p_{22}(X) &=
x_1^2+x_2^2-x_3^2-x_4^2,\\
p_{31}(X) &= x_3^2-x_4^2-2ix_3x_4,& \quad  p_{32}(X) &=
2(-x_1x_3+x_2x_4+i(x_1x_4+x_2x_3)), \\ 
p_{13}(X) &= x_3^2-x_4^2+2ix_3x_4 & \quad p_{23}(X) &=
x_1x_3+x_2x_4+i(x_1x_4-x_2x_3),
\end{align*} 
\[
 p_{33}(X) = x_1^2-x_2^2-2ix_1x_2.
\]

Finally, the action of the determinant representation is given by
\begin{equation}\label{det}
X(\alpha) = x_1^2 + x_2^2 + x_3^2 + x_4^2.
\end{equation}

Now, we are ready to define the Brandt matrices. Let $I_1 = \mathcal
O,\ldots,I_{12}$  be representatives of all the distinct left
$\mathcal O$--ideal  classes. Let $\mathcal O_j :=
I_j^{-1}I_j$ 
be the right  order of $I_j$ and let $e_j$ denote the number of units
in $\mathcal O_j = \frac{1}{Nr_{B/F}(I_j)}\bar I_jI_j$.
In other words,  $e_j$ is the number of times the positive definite
quadratic  form $Nr_{B/F}(x)$, $x \in \mathcal O_j$ represents $1$. 
Then, for any totally positive element $\xi \in \mathfrak o_F$ (\ie
$\xi \in \mathfrak o_{F}^+$)  we set
\begin{equation}\label{Brandt}
(B(\xi))_{ij} := e_j^{-1}\sum_{
\alpha\in I_j^{-1}I_i \atop 
Nr_S(\alpha)= \xi}X(\iota_1\alpha)\cdot[{^tX_2(\iota_2 \alpha)}]
\end{equation}

where $\{\iota_{\tau}\}_{\tau = 1,2} = I$ denotes the
set of the two real embeddings of $F = \Q(\sqrt 5)$ into $\R$. 
The symbol $Nr_S(\alpha)$ stands for the scaled norm of $\alpha
\in I_j^{-1}I_i$. If $J$ is an ideal of $B$,
we denote by $Nr(J)_+$ a totally
positive generator  of the (principal) ideal $Nr_{B/F}(J)$. The
scaled norm  of an element $\beta \in J$ is the totally
positive element $Nr_S(\beta) := \frac{Nr_{B/F}(\beta)}{Nr(J)_+}$. 
The $(B(\xi))_{ij}$'s are $(3\cdot 12)\times (3\cdot 12)$ square $\C$--matrices and they are  
divided  into $144$ $(3\times 3)$--blocks.

As $\xi$ varies among the set of totally positive elements of $\mathfrak
o_F$, the entries of the Brandt matrix series
\begin{equation}\label{theta}
\Theta(z) = \sum_{\xi\in\mathfrak
o_{F}^+}B(\xi)\exp(2\pi\sqrt{-1}~\xi~z),\quad z \in \mathfrak H =
\text{upper half plane}
\end{equation}
represent holomorphic Hilbert cusp forms of weight ${\bf k}$ on $S_{{\bf
k}}(\mathfrak P,5)$ on which the Brandt matrices act. In particular, these
matrices satisfy identities similar to those fulfilled by the Hecke operators. 

It is worth to recall the
following elementary property of the totally positive elements of $F$.
For a proof of it, we refer to any introductory book in algebraic
number theory ($w$ a fundamental unit in $F$)
\[
\mathfrak o_{F}^+ = \{\xi = a + bw~\vert~a,b\in \Z, a > 0,~\text{and}~
-\frac{a}{w} < b < \frac{a}{w-1}\}.
\]
Using this property, we can order the elements of $\mathfrak o_{F}^+$
lexicographically (\ie first by $a$ and then by $b$). Hence, the
sum in \eqref{theta} is well defined. 

It follows from the theory of the ``basis problem'' (\cf~\cite{HPS}) that 
$(B(\xi))_{ij}$ give a representation of the Hecke operators on the
space of cusp forms on the quaternion algebra $B$. Namely,
they describe the action of $(G_B)_\infty$
on $L_{\bf k}(\C)$ via the chosen map $i: (G_B)_\infty \hookrightarrow \mathfrak
G(\C)^I$ as we 
explained in the last paragraph.     
The explicit expression of these matrices
is subordinated to the knowledge of the left ideal
classes $I_j$, the right orders $\mathcal O_j$ and the associated
$\theta$--series. This part will be described in the next
section.

\section{The description of the algorithm.}\label{6}


Let $F = \Q(\sqrt 5)$ and let $w = \frac{1+\sqrt 5}{2}$ be a fundamental
unit. Let $B$ be the quaternion algebra with center $F$ 
 as defined in Lemma~\ref{quat}. We recall the definition of the
 Eichler order $\mathcal O$ of level $5$ defined in section~\ref{3}
\[
\mathcal O = \mathfrak o_F[1,~X,~-\frac{1}{2}+\frac{3-w}{2}Y,
-\frac{w}{2}-\frac{w+1}{2}X+\frac{1}{2}Y+\frac{w}{2}XY].
\]

The order $\mathcal O$ contains 12 different left-ideals classes whereas
the algebra $B$ has 3 distinct isomorphism classes
of orders with the same level as $\mathcal O$. 

The method used for finding the left-ideal classes is similar to that applied for
computing the ideal class group of a number field. The properties
$h(F) = 1$ and $Nr_{F/\Q}(w) = -1$ guarantee that any ideal of $\mathfrak o_F$
can be generated by a totally positive element and that every ideal of
$B$ has a basis over $\mathfrak o_F$. 

For a fixed $\alpha \in B \setminus F$, the quadratic field extension $K = F(\alpha)$ 
is contained in $B$. The non-trivial field
automorphism $\sigma$ of $K$ that fixes $F$ may be thought of as the
conjugation in $B$: $\alpha^\sigma = \bar\alpha$. In this way, the relative norm $Nr_{K/F}$
becomes the restriction $Nr_{\vert_{\bf k}}$ to $K$ of the norm $Nr_{B/F}$
defined on $B$. 

If $\mathcal I$ is an ideal of $\mathfrak o_K$ (the ring of integers of
$K$), the modulus $I = \mathcal O\mathcal I$ represents a left-ideal of $\mathcal O$. This is because
$\mathcal O(\mathcal O\mathcal I) = \mathcal O\mathcal
I$. Furthermore: $Nr_{B/F}(I) = Nr_{\vert_K}(\mathcal I)$, as $1 \in \mathcal
O$.  Finally, the representatives of the same class in the ideal class
group of $K$ give rise to elements that belong to the same class as left
ideals in $\mathcal O$. 

To an ideal in the quaternion algebra $B$, we associate  its
$\theta$--series. This series catalogues, in a way that we will  make clear
quite soon, the number of elements in the ideal with a fixed
norm. More precisely, if $I$ is a $\mathfrak o_F$--ideal of $B$, one
defines the $\theta$--series of $I$ as
\begin{align}\label{repn}
\theta_I(\tau) = &\sum_{\alpha\in I}\exp(\tau Nr_S(\alpha)) = \sum_{\xi
\in \mathfrak o_{F}^+}c_{\xi,I}\exp(\tau\xi),~\text{where}\\ &c_{\xi,I}
= \sharp\{\alpha \in I~\vert~Nr_S(\alpha) = \xi\}\notag. 
\end{align}

For the notations used here, we refer to \eqref{Brandt}. The symbols
$c_{\xi,I}$ are the representation numbers of $\xi$ in
$I$. Note that $c_{\xi,I}$ is finite for every $\xi$ and $I$. This is
an easy consequence of the fact that $B$ is totally
definite. For a fixed choice of $\xi$,
the definition of the representation numbers 
is independent of the choice of the positive generator of $Nr_{B/F}(I)$ (\ie
$Nr_+(I)$) involved in the definition of the scaled norm $Nr_S$
(\ie any two choices for $Nr_+(I)$ differ by a totally positive
unit in $F$). Thus, the $\theta$--series are well defined objects. 
These series are used to provide a necessary condition for testing
whether two
ideals in $B$ belong  to the same class. The proof of the following proposition is immediate.

\begin{prop} If $I$ and $J$ are $\mathfrak o_F$--ideals of $B$,
\[
J = \gamma_1 I\gamma_2,~\gamma_i \in B^\times\quad\text{then}\quad c_{\xi,J} =
c_{\xi,I}~\forall \xi \in \mathfrak o_{F}^+.
\]
\end{prop}

Notice that two ideals (or orders) in $B$ may have the
same $\theta$--series and generate different  classes. In order
to check whether two left $\mathcal 
O$-ideals belong to different  classes, one is sometimes required
to use the following necessary and sufficient condition.

\begin{prop}\label{iff} Let $I$ and $J$ be two left $\mathcal O$-ideals for an
Eichler order $\mathcal O$. Then, $I$ and $J$ belong to the same ideal
class if and only if there exists an element $\alpha \in \bar J I$
such that $Nr_{B/F}(\alpha) = Nr_{B/F}(I)Nr_{B/F}(J)$.
\end{prop}
\begin{proof}  We refer to \cite{P} Proposition~1.8, for the proof. 
The computations written down in \op generalize easily to any
quaternion algebra over a number field.\end{proof}

The criterium of Proposition~\ref{iff} is the most
convenient  for testing 
ideals and orders in $B$. It is clear though that its accuracy is
payed by time consuming! 

The following, is a concise description of the algorithm that we used in the process
 of testing ideal classes. This procedure requires to find the
representation numbers $c_{\xi}$ in the most
efficient way possible. The $c_{\xi}$'s describe the number of times the
quadratic form $Nr_S$ represents the totally positive definite element
$\xi \in \mathfrak o_F$. 

If $I$ is an ideal of $B$ with basis
$\{X_1,\ldots X_4\}$, any element $b \in I$ can be written in the form $b =
(b_{1}+b_{5}w)X_1 + (b_{2} + b_{6}w)X_2 + (b_{3}+b_{7}w)X_3 +
(b_{4}+b_{8}w)X_4$. Its scaled norm is then $Nr_S(b) = Q_1(Y)
+ Q_2(Y)w$, where $Y = [b_1,\ldots b_8] \in \Z^8$. The functions $Q_1$ and
$Q_2$ are quadratic forms with coefficients in $\Q$ and values in
$\Z$. Furthermore, $Q_1$ is
 (symmetric) positive definite because the quaternion algebra is totally definite. 
This property holds regardeless of
the chosen basis for $I$. 

If $\xi = \xi_1 + \xi_2w \in \mathfrak o_F$, then 
the finite set $\{Y \in \Z^8~\vert~Q_1(Y) = \xi_1\}$ (as
well as the similar set associated to $Q_2(Y)$) is
computable through a process of diagonalization of $Q_1(Y)$ obtained by
using the ``completion of square'' method. The efficiency
of this algorithm depends upon the choice of the basis for $I$. It can
be easily seen (and quite natural to expect) that the optimal choice for the basis
is the one associated to a Hermite normal form of the $4\times 4$
matrix that describes a basis for the ideal. The existence of such a form is
a consequence of the fact that $F$ is an Euclidean
domain (although the process is still possible in the more general
setting of a UFD domain). 

Said that, it is immediate to verify
that an efficient use of the necessary and sufficient condition of
proposition~\ref{iff} relies on the process of  first multiplying two
$4\times 4$ matrices (previously reduced in Hermite normal form) and
then using the Hermite normal form reduction to reduce a $4\times
16$ matrix into a $4\times 4$ upper--triangular. 

The main
part of the process is the search of the complete set of solutions for the
system: $Q_1(Y) = \xi_1 M_1$, $Q_2(Y) = \xi_2 M_2$, in the way 
explained before and for $Nr_{B/F}(I)Nr_{B/F}(J) = M_1 + M_2w$. This is
the most lengthy and tedious part. Its efficiency depends on a well-written
 implementable computer algorithm for finding the solutions of a
given quadratic form previously written as a sum of squares with positive integer
coefficients.

Now, we pass to the description of the left-ideals and their
right orders. 

We found the twelve left-ideals of $B$ by considering 
five distinct quadratic extensions $K = F(\alpha)$, $\alpha \in B$. We
do not claim  that the number and the choice of the extensions was
the best possible. It is quite plausible that one may find the whole bunch of classes by
considering a single extension $K$ as above.  
In the table below we have described these extensions together with   
their minimal polynomials $p_i$ over $F$ and the class numbers of
the quartic  extensions $K_i/\Q$. We denoted by $f_i$ $i = 1,\ldots 4$
the basis of the order $\mathcal O$ as defined in \eqref{basis}.
\begin{eqnarray*}
\begin{tabular}{|c|c|c|c|} 
$K_i$ & $\alpha_i \in \mathcal O$ & $p_i(x)$ & $h(K_i)$\\ 
\hline\hline 
$F(\alpha_1)$ & $f_2$ & $x^2+6$ & $4$ \\
\hline
$F(\alpha_2)$ & $f_2+3f_3-f_4$ & $x^2-5wx+85-18w$ & $22$\\  
\hline
$F(\alpha_3)$ & $(2-w)f_1-f_2-(w+1)f_3+wf_4$ & $x^2-(3-5w)x+17$ & $5$\\ 
\hline 
$F(\alpha_4)$ & $-f_2-(w+2)f_3+(2w-1)f_4$ & $x^2+5wx+31$ & $9$\\ 
\hline 
$F(\alpha_5)$ & $-(5w+2)f_1+2(2w+1)f_3+(1-4w)f_4$ & $x^2+wx+47$ & $17$  \\ 
\hline 
\end{tabular}
\end{eqnarray*}
We write
\begin{gather*}
\alpha_1 = X,\quad \alpha_2 =
\frac{5}{2}w+(\frac{1}{2}+w)X+\frac{5}{2}Y+\frac{1}{2}(3-w)XY, \\
\alpha_3 =
\frac{1}{2}(3-5w)-X-\frac{1}{2}(2+w)Y,\quad \alpha_4 =
-\frac{5}{2}w-X-\frac{5}{2}Y,\\
\alpha_5 =
-\frac{w}{2}+\frac{1}{2}(1-w)X+(\frac{5}{2}+2w)Y+\frac{1}{2}(-2+w)XY.
\end{gather*}

In the following table we have collected together the definitions of the
ideals $\mathcal I_i = (a_i+wb_i,\gamma_i)$ ($i= 1,\ldots,12$) in the 
quartic extensions $K_i$. We used the necessary and sufficient condition of
proposition~\ref{iff}  to
verify that  the associated left-$\mathcal O$-ideals
$\mathcal O\mathcal I_i = I_i$ define twelve distinct 
classes.  

\begin{eqnarray*}
\begin{tabular}{|c|c|c|c|c|} 
$\mathcal I_i$ & $K_i$ & $a_i+wb_i$ & $\gamma_i$ & $I_i \vert p \in
\Z$\\ \hline\hline 
$\mathcal I_1$ & $F$ & $1$ & & \\ \hline 
$\mathcal I_2$ & $F(\alpha_1)$ &
$2+w$ & $2+\alpha_1$ & $5$\\ \hline
$\mathcal I_3$ &
$F(\alpha_1)$ & $2$ & $2w+\alpha_1$ & $2$\\ \hline
$\mathcal I_4$ & $F(\alpha_1)$ & $2+w$ &
$3+\alpha_1$ & $5$\\ \hline 
$\mathcal I_5$ & $F(\alpha_2)$ & $2+w$
& $-2+\alpha_2$ & 
$5$\\ \hline   
$\mathcal I_6$ & $F(\alpha_2)$ & $7-2w$ & $16+\alpha_2$
 & $31$ \\ \hline  
$\mathcal I_{7}$ & $F(\alpha_2)$ & $9+13w$ &
$10+\alpha_2$& $29$\\ \hline 
$\mathcal I_8$ & $F(\alpha_3)$ & $4-w$ &
$8+\alpha_3$ & $11$\\ \hline
$\mathcal I_{9}$ & $F(\alpha_3)$ & $4-w$ &
$9+\alpha_3$ & $11$\\ \hline 
$\mathcal I_{10}$ & $F(\alpha_4)$ & $7$ &
$6+4w+\alpha_4$ & $7$\\ \hline 
$\mathcal I_{11}$ & $F(\alpha_5)$ & $3w+5$ &
$13+\alpha_5$ 
& $31$\\\hline 
$\mathcal I_{12}$ & $F(\alpha_5)$ & $17$ &
$12(1+w)+\alpha_5$ & $17$\\ \hline 
\end{tabular}
\end{eqnarray*}

The initial coefficients of the theta series of the ideals $I_i$
are tabulated below. Notice that since $\mathcal O$ has $2$ elements of 
norm $1$, every $\theta$--series coefficient  is a multiple of $2$.
  
\begin{eqnarray*}
\begin{tabular}{|c||c|c|c|c|c|c|c|c|c|c|c|c|}
\multicolumn{13}{c}{$\zeta = a + bw$}\\ \hline\hline
Ideals & $1$ & $2$ & $2+w$ & $3-w$ & $3$ & $3+w$ & $4-w$ & $4$ & $4+w$
& $4+2w$ & $5$ & $10$\\ \hline\hline 
$I_1 = \mathcal O$ & $2$ & $0$ & $0$ & $0$ & $0$ &
$4$ & $4$ & $2$ & $0$ & $0$ & $10$ & $20$ \\ \hline 
$I_2$ & $0$ & $0$ & $10$ & $10$ & $0$ & $2$ &
$2$ & $0$ & $4$ & $2$ & $18$ & $0$ \\ \hline 
$I_3$ & $0$ & $2$ & $0$ & $0$ & $2$ & $0$ &
$0$ & $0$ & $0$ & $0$ & $20$ & $10$ \\ \hline 
$I_4$ & $0$ & $0$ & $2$ & $2$ & $0$ & $2$ &
$2$ & $0$ & $4$ & $10$ & $0$ & $18$ \\ \hline 
$I_5$ & $0$ & $0$ & $2$ & $2$ & $0$ & $2$
& $2$ & $0$ & $4$ & $0$ & $10$ & $18$ \\ \hline 
$I_6$ & $0$ & $0$ & $0$ & $0$ & $0$ & $2$ &
$2$ & $0$ & $4$ & $2$ & $18$ & $10$  \\ \hline 
$I_7$ & $0$ & $0$ & $0$ & $0$ & $0$ & $2$ &
$2$ & $0$ & $4$ & $2$ & $18$ & $10$ \\ \hline 
$I_8$ & $0$ & $0$ & $0$ & $0$ & $0$ & $2$ &
$2$ & $0$ & $4$ & $2$ & $18$ & $10$ \\ \hline 
$I_9$ & $0$ & $0$ & $2$ & $2$ & $0$ & $2$ &
$2$ & $0$ & $4$ & $0$ & $10$ & $18$ \\ \hline 
$I_{10}$ & $0$ & $0$ & $2$ & $2$ & $0$ & $2$ &
$2$ & $0$ & $4$ & $0$ & $10$ & $18$ \\ \hline 
$I_{11}$ & $0$ & $0$ & $0$ & $0$ & $0$ & $2$
& $2$ & $0$ & $4$ & $2$ & $18$ & $10$ \\ \hline 
$I_{12}$ & $0$ & $0$ & $2$ & $2$ & $0$ & $2$
& $2$ & $0$ & $4$ & $0$ & $10$ & $18$ \\ \hline 
\end{tabular}
\end{eqnarray*}

The following table shows few of the coefficients of the $\theta$--series  of the
right orders $\mathcal O_i := I_i^{-1}I_i$. 
\begin{eqnarray*}
\begin{tabular}{|c||c|c|c|c|c|c|c|c|c|c|c|c|c|c|c|c|}
\multicolumn{15}{c}{$\zeta = a + bw$}\\ \hline\hline
Right-Orders & $e_i$ & $2$ & $3$ & $4$ & $5$ & $6$ & $7$ & $8$
& $9$ & $10$ & $11$ & $11+w$ & $12-w$ & $13+w$\\ \hline\hline 
$\mathcal O_1 = \mathcal O$
& $2$ & $0$ & 
$0$ & $2$ & $10$ & $2$ & $20$ & $0$ & $2$ & $20$ & $28$ & $4$ & $4$ &
$44$ \\ \hline  
$\mathcal O_2$ &  $2$ & $0$ & $0$
& $2$ & $10$ & $2$ & $20$ & $0$ & $2$ & $20$ & $28$ & $4$ & $4$ & $44$ \\ \hline 
$\mathcal O_3$ &  $2$ & $0$ & $0$
& $2$ & $10$ & $2$ & $20$ & $0$ & $2$ & $20$  & $28$ & $4$ & $4$ & $44$
\\ \hline  
$\mathcal O_4$ &  $2$ & $0$ & $0$
& $2$ & $10$ & $2$ & $20$ & $0$ & $2$ & $20$ & $28$ & $4$ & $4$ & $44$ \\ \hline 
$\mathcal O_5$ &  $2$ & $0$ & $0$ &
$2$ & $26$ & $0$ & $8$ & $0$ & $2$ & $0$ & $20$ & $14$ & $16$ & $20$\\ \hline 
$\mathcal O_6$ &  $2$ & $0$ & $0$ &
$2$ & $26$ & $0$ & $8$ & $0$ & $2$ & $0$ & $20$ & $16$ & $14$ & $16$\\ \hline 
$\mathcal O_7$ &  $2$ & $0$ & $0$ &
$2$ & $26$ & $0$ & $8$ & $0$ & $2$ & $0$ & $20$ & $16$ & $14$ & $16$\\ \hline 
$\mathcal O_8$ &  $2$ & $0$ & $0$
& $2$ & $26$ & $0$ & $8$ & $0$ & $2$ & $0$ & $20$ & $16$ & $14$ & $16$ \\ \hline 
$\mathcal O_9$ &  $2$ & $0$ & $0$
& $2$ & $26$ & $0$ & $8$ & $0$ & $2$ & $0$ & $20$ & $14$ & $16$ & $20$ \\ \hline 
$\mathcal O_{10}$ &  $2$ & $0$ & $0$ &
$2$ & $26$ & $0$ & $8$ & $0$ & $2$ & $0$ & $20$ & $16$ & $14$ & $16$ \\ \hline 
$\mathcal O_{11}$ &  $2$ & $0$
& $0$ & $2$ & $26$ & $0$ & $8$ &  $0$ & $2$ & $0$ &  $20$ & $14$ & $16$
& $20$ \\ \hline 
$\mathcal O_{12}$ &  $2$ & $0$
& $0$ & $2$ & $26$ & $0$ & $8$ & $0$ & $2$ & $0$ & $20$ & $14$ & $16$ & $20$ \\ \hline 
\end{tabular}
\end{eqnarray*} 

Because $t(\mathcal O) = 3$,  there are only three
different types (\ie conjugacy classes by elements of $B^\times$) of
orders of $B$ with level 5. One
can classify them by 
looking at the corresponding $\theta$--coefficients. The above table 
 shows that  at $\xi = 5$ one can distinguish the first
equivalence class. At $\xi = 11 + w$ the three classes become
all distinct. We  have chosen $\mathcal O$,
$\mathcal O_5$ and $\mathcal O_6$ as representatives for these classes. 
The lattices spanned by $\bar{I_5}I_5
$ and by $\bar{I_6}I_6$ are given by the following upper triangular
matrices (each column describes the coefficients of \resp
$1,~X,~Y,~XY$). The ratio in front of each matrix represents the norm
of the related $\mathcal O$--left--ideal (notice that $\mathcal O_j
= \frac{1}{Nr_{B/F}(I_j)}\bar I_j I_j$)  
\begin{multline*}
\mathcal O_5 = \mathfrak o_F\left [\frac{1}{2+w}\begin{array}({cccc})
4w+3 & 0 & -(2w+\frac{3}{2}) & -\frac{1}{2}(11w+7) 
\\ 0 & 5w & -(3w+1) & w-\frac{1}{2} \\
0 & 0 & -\frac{5}{2}(w+1) & -(\frac{3}{2}w+1) \\
0 & 0 & 0 & \frac{1}{2}(w-1) 
\end{array}\right ],\\\\ 
\mathcal O_6 = \mathfrak o_F\left [\frac{1}{7-2w}\begin{array}({cccc})
2w-7 & 0 & -w+\frac{7}{2} & 0 
\\ 0 & 73w+39 & 30w+19 & -(\frac{55}{2}w+19) \\
0 & 0 & \frac{1}{2}(11w-23) & 2(w-2) \\
0 & 0 & 0 & w+\frac{1}{2} 
\end{array}\right ].
\end{multline*}

For the computation of the Brandt matrices $(B(\xi))_{i,j}$ 
we needed to know, for the first few  
totally positive $\xi \in \mathfrak o_{F}$, the number of elements
$\alpha \in I_j^{-1}I_i$ with scaled norm $\xi$. Due to the symmetry
of the Brandt matrices (\cf~\cite{P}), it suffices to compute
these numbers for the 78 ideals
$I_j^{-1}I_i$ for $j \ge i$. 

Because the Brandt matrices  
generate a commutative semisimple ring it is possible to diagonalize
them simultaneusly. This means that one can choose a finite set of
totally positive elements $\xi$ and for each of them compute
the corresponding Brandt matrices and then search
for a simultaneus diagonalization. In our case we only needed to verify that among
the eigenvectors there was one whose eigenvalues mapped down to $\Q$ via the norm map,
match with the traces of the Frobenius listed in paragraph~\ref{J2}.

Note also that it is enough to compute 
the characteristic polynomial (of degree 36) of the matrix 
$^t{X_2(\iota_2(\alpha))}$ (\cf~\eqref{Brandt}) for a particular choice of $\xi$ and then 
carry on the computations for it. 
We decided to choose
$\xi = 3+w$ (one of the two primes in
$F$ above 11 ) and factor the related polynomial in $F$.  
This polynomial splits in 10 linear factors (each
of which appears with multiplicity greater than one)
and a factor of degree 3 (counted with multiplicity 2). 
Note that one has to take into account the presence of the determinant representation in
the definition of  $(B(\xi))_{ij}$. In fact, it turned out that
one of the
$\mathfrak o_F$--roots of the associated Brandt polynomial (appearing
with  multiplicity $3$) is
$-4(4+w)(4-w) = 4(-15+w)$. The factor $-4(4+w)$ represents the
root associated to the representation $X_2$, whereas $4-w$ describes
the action of the determinant. The first number coincides with the 
degree-one coefficient of one of the two quadratic factors in which  the
characteristic polynomial of the Frobenius at 11 splits (\cf section~\ref{J2}).

Using the same technique we computed the characteristic polynomial of
the matrix $(B(\xi))_{i,j}$ at the inert prime $\xi = 7$. This polynomial has
16 rational roots (counted with multiplicities). One of these is
$-10\cdot 7 = -70$ (with multiplicity 4). Here, $-10$ is the root
associated to $X_2$. Again, the square of this number
coincides with the degree-two coefficient in the
characteristic polynomial of the Frobenius at 7.

At this point, we considered the common eigenspace of the two matrices at the
eigenvalue of interest. This is a one-dimensional space.
Its eigenvector $\varepsilon$ was the natural candidate for
representing the Hilbert modular form $\mathfrak f$ (\cf paragraph~\ref{4}). To compute
further eigenvalues associated to $\varepsilon$, it was enough to work
with a particular column of the Brandt matrix and multiply that with 
$\varepsilon$. 

The following is the description of the (transpose of the)
eigenvector $\varepsilon$ (we denote by $v = \sqrt{3-w}$): \newline\newline
$
\varepsilon = [
0,-w - 2,0,
0,1,0,
0,-\frac{1}{2}w-1,0,
0,1,0,
(-\frac{3}{25}w+\frac{4}{25})v\sqrt{-6}+(\frac{3}{50}w-\frac{2}{25})v,\frac{1}{1
0}w-\frac{3}{10},(-\frac{3}{25}w+\frac{4}{25})v\sqrt{-6}+(-\frac{3}{50}w+\frac{2
}{25})v,
(\frac{833}{4805}w+\frac{486}{4805})v\sqrt{-6}+(\frac{467}{9610}w +
\frac{142}{4805})v,-\frac{183}{1922}w - \frac{101}{1922},(\frac{833}{4805}w +
\frac{486}{4805})v\sqrt{-6} + (-\frac{467}{9610}w - \frac{142}{4805})v,
(\frac{63}{4205}w -\frac{4}{4205})v\sqrt{-6} + (\frac{847}{8410}w -
\frac{521}{8410})v,-\frac{2}{841}w + \frac{67}{1682},(\frac{63}{4205}w -
\frac{4}{4205})v\sqrt{-6} + (-\frac{847}{8410}w + \frac{521}{8410})v,
(\frac{4}{605}w +\frac{23}{605})v\sqrt{-6} + (-\frac{169}{1210}w +
\frac{87}{1210})v,-\frac{2}{121}w - \frac{23}{242},(\frac{4}{605}w
+\frac{23}{605})v\sqrt{-6} + (\frac{169}{1210}w - \frac{87}{1210})v,
(\frac{89}{605}w +\frac{58}{605})v\sqrt{-6} + (\frac{441}{1210}w +
\frac{267}{1210})v,\frac{60}{121}w + \frac{85}{242},(\frac{89}{605}w +
\frac{58}{605})v\sqrt{-6} + (-\frac{441}{1210}w - \frac{267}{1210})v,
(\frac{4}{245}w + \frac{3}{245})v\sqrt{-6} + (\frac{111}{490}w +
\frac{31}{245})v,-\frac{5}{98}w - \frac{5}{98},(\frac{4}{245}w +
\frac{3}{245})v\sqrt{-6} + (-\frac{111}{490}w - \frac{31}{245})v,
(\frac{19}{4805}w + \frac{63}{4805})v\sqrt{-6} + (\frac{1221}{9610}w -
\frac{884}{4805})v,\frac{107}{1922}w - \frac{151}{1922},(\frac{19}{4805}w +
\frac{63}{4805})v\sqrt{-6} + (-\frac{1221}{9610}w + \frac{884}{4805})v,
(-\frac{1}{85}w + \frac{16}{1445})v\sqrt{-6} + (-\frac{33}{2890}w -
\frac{11}{2890})v,\frac{4}{289}w - \frac{19}{578},(-\frac{1}{85}w
+\frac{16}{1445})v\sqrt{-6} + (\frac{33}{2890}w + \frac{11}{2890})v].
$\vspace{.2in}

The following table resumes the first few eigenvalues of
$\varepsilon$ when $\xi$ is a prime (one can compare these values with those listed
in table $\eqref{tracetable}$).
\begin{eqnarray*}
\begin{tabular}{|c|c|c|} 
$\xi = a + wb$ & $\xi \vert p \in \Z$ & $\varepsilon$\\ \hline\hline 
$1$ & $$ & $1$ \\ \hline 
$2$ & $2$ &
$4$\\ \hline
$3$ &
$3$ & $9$ \\ \hline
$-1+2w$ & $5$ & $0$
\\ \hline   
$7$ & $7$ & $-70$ 
\\ \hline 
$3+w$ & $11$ & $4(-15+w)$\\
$4-w$ & $11$ & $-4(14+w)$ 
\\\hline 
$13$ & $13$ & $2990$ 
\\ \hline 
$17$ & $17$ & $-170$ 
\\ \hline 
$4+w$ & $19$ & $4(12-29w)$\\
$5-w$ & $19$ & $-4(17-29w)$ 
\\ \hline 
$23$ & $23$ & $3450$ 
\\ \hline 
$5+w$ & $29$ & $2(19-8w)$\\
$6-w$ & $29$ & $2(11+8w)$ 
\\ \hline 
$7-2w$ & $31$ & $24(3-5w)$ \\
$5+2w$ & $31$ & $-24(2-5w)$
\\ \hline 
$37$ & $37$ & $-29970$ 
\\ \hline 
$7-w$ & $41$ & $-2(13+132w)$\\
$6+w$ & $41$ & $-2(145+132w)$ 
\\ \hline 
$43$ & $43$ & $149210$ 
\\ \hline 
$47$ & $47$ & $93530$ 
\\ \hline 
$53$ & $53$ & $-235850$ 
\\ \hline 
$2w+7$ & $59$ & $-4(141+8w)$\\
$9-2w$ & $59$ & $4(149-8w)$\\ \hline 
\end{tabular}
\end{eqnarray*}

\begin{thm}
For each $\xi$ the eigenvalue of $B(\xi)$ corresponding to the above
eigenvector is in $F$.
\end{thm}
\begin{proof}
Suppose 
$\alpha=x_1 + x_2\cdot{\bf I}
+ x_3\cdot{\bf J}+ x_4\cdot{\bf K}$ 
is an element of $B$. Then $x_1\in F$, $x_2\in F\cdot\sqrt{6}$, 
$x_3\in F\cdot v$ and $x_4\in F\cdot v\sqrt{6}$.
Recall the definition of $X_2(\alpha)$ from $\eqref{sym}$ in terms
of harmonic polynomials $p_{ij}$. From the definition of
the $p_{ij}$ it immediately follows that 
$X_2(\alpha)_{ij}\in F(\sqrt{-6})$ if $(i,j)=(1,1)$, $(1,3)$,
$(3,1)$, $(3,3)$ or $(2,2)$, and $X_2(\alpha)_{ij}\in
F(\sqrt{-6})\cdot v$ otherwise. From the definition of
$B(\xi)$ it follows that $B(\xi)_{ij}\in F(\sqrt{-6})$ if
$(i,j)$ is congruent to $(1,1)$, $(1,3)$,
$(3,1)$, $(3,3)$ or $(2,2)$ mod 3, and 
$B(\xi)_{ij}\in F(\sqrt{-6})\cdot v$ otherwise.
Note that the eigenvector $\varepsilon$ satisfies
$\varepsilon_i\in F(\sqrt{-6})$ if $i\equiv 2\ \mod\ 3$, and
$\varepsilon_i\in F(\sqrt{-6})\cdot v$ otherwise. So 
$B(\xi)\varepsilon$ also has this property, hence the eigenvalue
of $\varepsilon$ is in $F(\sqrt{-6})$. 

It is known that eigenvalues of the Hecke action on Hilbert modular
forms are contained in a totally real field. This implies that
the eigenvalue of $\varepsilon$ is in $F$.
\end{proof}

\ifx\undefined\bysame
\newcommand{\bysame}{\leavevmode\hbox to3em{\hrulefill}\,}
\fi

\end{document}